\documentclass[twocolumn]{autart}    % Enable this line and disable the 
\usepackage{amsfonts}
\usepackage[all]{xy}
\usepackage{amsmath}
\usepackage{mathrsfs,amssymb,amsmath,booktabs,array,xcolor,url,color,soul,multirow}
\usepackage{mathbbold,bm}
\usepackage{float}
\usepackage{graphicx}
\usepackage{arydshln}
\usepackage{subfigure}
\usepackage{enumitem}
\usepackage{ntheorem}

\usepackage{soul}
\soulregister\cite7
\soulregister\ref7 
\usepackage{wrapfig}

\usepackage{graphicx}
\usepackage{wrapfig}
\usepackage[round]{natbib}

\newtheorem{theorem}{Theorem}[section]

\newtheorem{lemma}[theorem]{Lemma}

\newtheorem{example}{Example}
\newtheorem{definition}[theorem]{Definition}
\newtheorem{remark}{Remark}     
\begin{document}
\begin{frontmatter}

\title{Network Combination to Persistence of High-dimensional Delayed Complex Balanced Mass-action Systems\thanksref{footnoteinfo}} % Title, preferably not more % than 10 words.
\thanks[footnoteinfo]{ This
article was presented in part at the paper accepted by 62nd IEEE Conference on Decision and Control, Dec. 13-15, 2023, Singapore. Corresponding author C. Gao. Tel. +86-571-87952431. 
Fax +86-571-87953794.}
\author[1]{Xiaoyu Zhang}\ead{Xiaoyu\_Z@zju.edu.cn},    % Add the 
%\author[1]{Shibo He}\ead{s18he@zju.edu.cn},         % e-mail address
\author[2]{Chuanhou Gao}\ead{gaochou@zju.edu.cn},  % (ead) as shown
\author[4]{Denis Dochain}\ead{denis.dochain@uclouvain.be}
\address[1]{Department of Control Science and Engineering, Zhejiang
University, Hangzhou 310027, China}
%\address[3]{Department of Biosystems Science and Engineering, ETH Zurich, 8092 Zurich, Switzerland}
\address[2]{School of Mathematical Sciences, Zhejiang
University, Hangzhou 310027, China}  % Please supply 
\address[4]{ICTEAM, UCLouvain, B\^{a}timent Euler, avenue Georges lema\^{i}tre 4-6, 1348 Louvain-la-Neuve, Belgium} 
           % full addresses
       % here.
                             % keyword list or with the 

\begin{abstract}
Complex balanced mass-action systems (CBMASs) are of great importance in the filed of biochemical reaction networks. However analyzing the persistence of these networks with high dimensions and time delays poses significant challenges. To tackle this, we propose a novel approach that combines 1-dimensional (1d) or 2d delayed CBMASs (DeCBMASs) and introduces inheritable combination methods based on the relationship between semilocking sets and intersecting species. These methods account for various scenarios, including cases where the set of intersecting species is empty, or there are no common species in non-trivial semilocking sets and the intersecting species set, or when special forms are present.
By utilizing these combination methods, we derive sufficient conditions for the persistence of high-dimensional DeCBMASs. This significantly expands the known class of delayed chemical reaction network systems that exhibit persistence. 
The effectiveness of our proposed approach is also demonstrated through several examples, highlighting its practical applicability in real-world scenarios. This research contributes to advancing the understanding of high-dimensional DeCBMASs and offers insights into their persistent behavior.
%and then design a controller
%that meets the expectations. }
\end{abstract}
\end{frontmatter}
\section{Introduction} 
Persistence is a vital performance indicator that reflects the coexistence of species within a system. Initially proposed in the context of ecosystems, it is closely tied to biodiversity \citep{Bertuzzo2011, Fischer2019}. In the field of biochemistry reaction systems, persistence plays a central role in determining the system ability to maintain functionality over time despite the presence of complex biochemical interactions. As such, persistence holds deep biological significance and has drawn considerable attention. In the theoretical study of biochemical reaction systems, persistence can be interpreted as a property that prevents trajectories from converging to boundaries. This can be evaluated by examining the presence of boundary points within the $\omega$-limit set of each bounded trajectory. However research on the persistence of such networks continues to face the double challenges of nonlinearity and high dimensionality. Moreover the inherent uncertainty associated with system parameters poses a significant obstacle to perform accurate analysis. The complexity of high-dimensional nonlinear systems often gives rise to intricate dynamical behaviors such as multistability and limit cycles, where even slight perturbations can lead to drastic changes in system dynamics. Thus studying robust persistence with a focus on network structure represents an important and long-term objective in the field of systems biology. 

Centered around this objective, a well-known conjecture called ``Persistence Conjecture," proposed by \cite{Feinberg1987}, asserts that all networks with weakly reversible structures are persistent. Drawing inspiration from this conjecture, significant efforts and progress have been made in the field. \cite{Angeli2007} and \cite{Anderson2008} introduced a specialized subset of species sets, referred to as siphon and semilocking sets (which are essentially equivalent). These sets share a notable characteristic: once a trajectory enters their respective boundaries (where the concentration of species in the semilocking set is zero), it remains indefinitely within them. Furthermore they demonstrated that the $\omega$-limit point of each trajectory cannot exist in the non-semilocking set. The introduction of semilocking sets in the analysis of biochemical reaction networks marks a significant milestone with deep implications for the study of persistence. It has also paved the way for further exploration by highlighting a promising direction: investigating the attractability of points situated on the boundaries associated with semilocking sets.
Building upon this foundation, numerous remarkable findings have emerged. These include the examination of complex balanced networks with low-dimensional stoichiometric subspaces \citep{Craciun2013}, or with specific types of semilocking sets \citep{Anderson2008, Anderson2010}, or with a single linkage class \citep{Anderson2011}; the utilization of conservation relationships among species within semilocking sets \citep{Angeli2007}; as well as investigations into low-dimensional endotactic networks \citep{Craciun2013, P2012} 
 and strong endotactic networks \citep{Anderson2020, G2014}.

Apart from the uncertainties in parameters, biological systems also face unpredictability in complex intermediate processes (some intermediate processes in biological systems are intricate, transient, and challenging to observe). The substantial complexity of intermediate processes also results in a major increase of system dimensionality. Consequently traditional modeling approaches encounter limitations in analysis and become less reliable in this case. Therefore when intermediate species are not the focus of interest, it is possible to simplify the modeling and reduce dimensionality significantly by introducing time delays to replace all intermediate processes like typical enzyme catalysis processes or DNA expression processes. In addition, time delays are essential for the generation and maintenance of oscillatory behavior in biochemical systems \citep{B2005}. They also play a crucial role in stabilizing the system dynamics \citep{F2014}. Therefore the study of persistence in biochemical reaction network systems with time delays carries important practical significance and application value. 
Hangos et al. \citep{G2018} were pioneers in recognizing the importance of extending and studying chemical reaction network theory with time delays. They introduced relevant definitions and formulated the Lyapunov functionals for complex balanced systems with delays by using chain method in \cite{Repin1965}. Their work revealed that such systems can maintain their local asymptotic stability even in the presence of time delays \citep{G2018b}. \cite{Zhang2022} further study the local asymptotic stability of two kinds of delayed chemical reaction systems relative to other invariant classes by conducting the decomposition of the solution space. Aside from the stability, the persistence of the delayed system also been researched. \cite{H2019} focuses on this topic and derives that the conservative delayed chemical reaction network with conservation relations between the species in semilocking set is persistent. \cite{Zhang2021} derives the persistence of the delayed complex balanced mass action system (DeCBMAS) with 2d stoichiometric subspace.

It is evident that the analysis of persistence based on network structure faces significant challenges, especially in high-dimensional systems where the stoichiometric subspaces of complex balanced systems exceed dimension 3. To address this issue, we propose novel inheritance combination methods based on the relationship between semilocking sets and the set of intersecting species denoted as $\mathcal{S}_c$. These methods involve integrating lower-dimensional sub-DeCBMASs, such as 1-dimensional or 2-dimensional systems, to gain insights into the persistence of high-dimensional systems.
When there are no intersecting species ($\mathcal{S}_c=\emptyset$), the coupling between subsystems is absent, making the persistence of high-dimensional DeCBMASs trivial. Therefore, our focus is on the case where $\mathcal{S}_c\neq \emptyset$, and we address the coupling between sub-systems.
First we consider the scenario where the set of intersecting species does not exist in any non-trivial semilocking set ($\mathcal{S}_c\cap W=\emptyset$). For the non-empty case, we further investigate several specific forms. By treating each semilocking set $W$ of the high-dimensional DeCBMASs as combinations of vertex-type and facet-type semilocking subsets, we derive sufficient conditions of persistence without relying on the parameters of high-dimensional DeCBMASs.
This research significantly expands the class of delayed chemical reaction systems with persistence. Furthermore we demonstrate the effectiveness of our proposed inheritance results through multiple examples. The derived results provide valuable insights into the dynamical analysis of high-dimensional systems, overcoming previous limitations.

This paper is organized as follows. Section \ref{sec:2} introduces the basic concepts about the DeCBMASs and persistence.
Section \ref{sec:3} gives the problem motivation of this paper. Section \ref{sec:4} focuses on the persistence of higher-dimensional DeCBMAS composed of 2d sub-DeCBMASs where the intersecting species does not exist in non-trivial semilocking set. Section \ref{sec:5} make some try in the higher-dimensional DeCBMASs where the intersecting species also exists in non-trivial semilocking sets. Section \ref{sec:6} concludes this paper. 
~\\~\\
\textbf{Notations:} $\bar{\mathscr{C}}_{+}=C([-\tau_{\mathrm{m}},0];\mathbb{R}^{n}_{\geq 0})$ is the space of non-negative function vectors defined on the interval $[-\tau_{\mathrm{m}},0]$, and $\mathscr{C}_+=C([-\tau_{\mathrm{m}},0];\mathbb{R}^{n}_{> 0})$ is the corresponding positive solution space; $\text{supp}~x$ denotes the support set of a vector $x$ defined as $\text{supp}~ x=\{S_j|x_{j}\neq 0\},~\forall x\in \mathbb{R}^{n}_{\geq 0}$; $\pi_{W}(v)$ denotes the projection of $v$ onto the components $v_j$ with $X_j\in W$.

\section{Preliminaries}\label{sec:2}
Prior to presenting our primary findings, we shall introduce some fundamental definitions related to chemical reaction networks (CRNs)\citep{Feinberg2019} and the delayed mass action kinetics \citep{G2018}, respectively.

%Let $\mathcal{N}$ be a chemical reaction network with $n$ species and $r$ reactions. Each reaction $R_i$ can be expressed in terms of species $X_j$ and their corresponding stoichiometric coefficients, as shown below:

Consider chemical reactions $R_i~ (i=1,...,r)$
\begin{equation} \label{eq:1}
R_i: \quad\sum^{n}_{j=1}y_{ji}S_j\stackrel{}{\longrightarrow} \sum^{n}_{j=1}y'_{ji}S_j,
\end{equation}
where $y_{ji}, y'_{ji}\in\mathbb{Z}^n_{\geq 0}$ represent the stoichiometric coefficients of species $S_j$. Denote by $y_{.i}=(y_{1i},...,y_{ni})^\top$ and $y'_{.i}=(y'_{1i},...,y'_{ni})^\top$ the reactant and product complexes, respectively. Then a CRN $\mathcal{N}$ can be uniquely determined by the triple $(\mathcal{S},\mathcal{C},\mathcal{R})$, where $\mathcal{S}$ represents the set of $n$ species, $\mathcal{C}$ denotes the set of all reactant and product complexes, and $\mathcal{R}$ represents the set of $r$ reactions, usually labelled by $\mathcal{N}=(\mathcal{S},\mathcal{C},\mathcal{R})$. 

Once the reaction $R_i$ occurs, the changes of all species can be described by the reaction vector defined as $y'_{.i}-y_{.i}$, which further forms a subspace of $\mathbb{R}^n$.
%can be defined, which reflects the structure of the network from a deeper perspective in terms of space.
\begin{definition}\em(stoichiometric subspace)
The stoichiometric subspace of $\mathcal{N}=(\mathcal{S},\mathcal{C},\mathcal{R})$ is defined as the linear span of the vectors $y'_{.i}-y_{.i},~ (i=1,\dots,r)$, given by
\begin{equation}
\mathscr{S} = \mathrm{span}\{y'_{.i}-y_{.i}|i=1,\dots,r\},
\end{equation}
whose dimension, denoted by $\dim \mathscr{S}$, is called the dimension of $\mathcal{N}$. The orthogonal complement of $\mathscr{S}$, denoted by $\mathscr{S}^\bot$, is given by
\begin{equation}
\mathscr{S}^\bot = \{a\in\mathbb{R}^n|a^\top y=0~\text{for all}~y\in\mathscr{S}\}.
\end{equation}
\end{definition}

%To model the evolution of a chemical reaction system, the rates of reactions, rather than just the network structure, must be considered.
It is commonly to use mass-action kinetics to evaluate reaction rate that follows a power law according to the concentration of each reactant species. 
\begin{definition}\em(mass-action law)
For each reaction $R_i$ in $\mathcal{N}$, it is said to obey mass-action law if its reaction rate $\delta_i(x)$ is evaluated by
 \begin{equation}
\delta_i(x)=k_i\prod^{n}_{j=1}x^{y_{ji}}\triangleq k_ix^{y_{.i}},
\end{equation}
where $x_j\in \mathbb{R}_{\geq 0}$ is the concentration of species $S_j(j=1,...,n)$, and the positive real number $k_i$ is called the reaction rate constant.
\end{definition}

An $\mathcal{N}=(\mathcal{S},\mathcal{C},\mathcal{R})$ with mass-action kinetics is called a mass-action system (MAS), labelled by $\mathcal{M}=(\mathcal{S},\mathcal{C},\mathcal{R},k)$ in the context. If there are reactions with time delays, meaning that reactant consumption is instantaneous while product generation is delayed, the MAS is named as the delayed mass-action system (DeMAS), denoted by $\mathcal{M}_{\mathrm{De}}=(\mathcal{S,C,R},k,\tau)$, where $k,\tau$ are the vectors of reaction rate constants and of time delays, respectively. The dynamics of $\mathcal{M}_{\mathrm{De}}$ can be represented by a set of delay differential equations
\begin{equation} \label{eq:dde}
\dot{x}(t)=\sum^{r}_{i=1}k_i[x(t-\tau_i)^{y_i}y'_{.i}-x(t)^{y_{.i}}y_{.i}],\quad t\geq 0
\end{equation}
where $x=(x_1,...,x_n)^\top$ is the vector of species concentrations, and $\tau_i\geq 0$ are constant time delays. Let $\tau_{\text{m}}$ be the maximum among all $\tau_i$, then the solution space of \eqref{eq:dde} is $\bar{\mathscr{C}}_{+}=C([-\tau_{\mathrm{m}},0];\mathbb{R}^{n}_{\geq 0})$. Based on the mass conservation relationships \citep{G2018}, $\bar{\mathscr{C}}_{+}$ can be decomposed into infinitely many equivalent classes.
\begin{definition}\em(stoichiometric compatibility class)\label{def:scc} For an $\mathcal{M}_{\mathrm{De}}=(\mathcal{S,C,R},k,\tau)$ with the dynamics of \eqref{eq:dde} and the solution space $\bar{\mathscr{C}}_+$, a non-negative (positive) stoichiometric compatibility class containing $\psi\in \bar{\mathscr{C}}_+$ is defined as:
\begin{small}
\begin{equation}\label{eq:scc}
\begin{split}
\mathcal{D}_\psi &= \{\theta\in \bar{\mathscr{C}}_+|a^{\top}g(\theta)=a^{\top}g(\psi)~{\rm for~all}~a\in\mathscr{S}^\bot\},\\
(\mathcal{D}^+_\psi &= \{\theta\in \mathscr{C}+|a^{\top}g(\theta)=a^{\top}g(\psi)~{\rm for~ all} ~a\in\mathscr{S}^\bot\},)
\end{split}
\end{equation}
where $g:\bar{\mathscr{C}}_+\rightarrow \mathbb{R}^{n}$ is defined by:
\begin{equation}\label{eq:zc0}
\begin{split}
g(\psi)&=\psi(0)+\sum^r_{i=1}\left(k_i\int^0_{-\tau_i}\psi(s)^{y_i}ds\right)y_i.\\
\end{split}
\end{equation}
\end{small}
\end{definition}
\cite{G2018b} proved that each trajectory $x^\psi(t)$ with an initial point $\psi$ will always stays in $\mathcal{D}_{\psi}$, namely, the stoichiometric compatibility class is an invariant set for each trajectory in $\mathcal{M}_{\mathrm{De}}$.

%The integral of the above equation means that the solution $x(t)$ with respect to some initial point $x_0\in\mathbb{R}^n_{\geq 0}$ will always stay in a invariant set, termed by the non-negative stoichiometric compatibility class following
%\begin{equation*}
%\bar{\mathcal{P}}^+_{x_0}=\{x(t)\in\mathbb{R}^n_{\geq 0}\;|\;x(t)-x_0\in \mathscr{S}\}.
%\end{equation*}

%Then the stoichiometric compatibility class $\bar{\mathcal{P}}^+_{x_{0}}$ can be rewritten as
%\begin{equation}\label{eq:wtd}
 % \bar{\mathcal{P}}^+_{x_{0}}=\{x\in\mathbb{R}^n_{\geq 0}\vert v^{\top}x=v^{\top}x_{0},\;\forall v\in \mathscr{S}^{\bot}\}.
%\end{equation}
\begin{definition}\em(complex balanced equilibrium)
A positive vector $\bar{x}\in \mathbb{R}^n_{>0}$ is called a \textbf{positive equilibrium} in $\mathcal{M}_{\mathrm{De}}=(\mathcal{S,C,R},k, \tau)$ if it satisfies $\dot{\bar{x}}=0$ in (\ref{eq:dde}). Furthermore, $\bar{x}$ is called a \textbf{complex balanced equilibrium} in $\mathcal{M}_{\mathrm{De}}$ if, for any complex $\eta\in\mathbb{Z}^n_{\geq 0}$ in the network, it satisfies the complex balance condition:
\begin{equation}\label{eq:cb}
\sum_{i:~y_{.i}=\eta}k_i\bar{x}^{y_{.i}}=\sum_{i:~y'_{.i}=\eta}k_i\bar{x}^{y_{.i}}.
\end{equation}
\end{definition}
A DeMAS possessing a complex balanced equilibrium is called a DeCBMAS. Equation \eqref{eq:cb} indicates that at any complex balanced equilibrium $\bar{x}$, the inflow and outflow rates of each complex $\eta$ are equal. Therefore the above definition demonstrates that the system $\mathcal{M}_{\mathrm{De}}$, at each complex balanced equilibrium $\bar{x}$, achieves a balance not only among the species but also among the complexes.

%A chemical reaction system having a complex balanced equilibrium is called a complex balanced system. \cite{Horn1972} proved that each non-negative stoichiometric compatibility class contains a unique complex balanced equilibrium for a 
%A positive vector $\bar{x}\in \mathbb{R}^n_{>0}$ is called a positive equilibrium if it satisfies $\sum^{r}_{k=1}k_k\bar{x}^{y_k}(y'_k-y_k)=0$ in Eq. (\ref{eq:2}); and $\bar{x}$ is called a \emph{complex balanced} equilibrium if for any complex $\eta\in\mathbb{Z}^n_{\geq 0}$ in the network it satisfies
% \begin{equation*}
%\sum_{k:~y_k=\eta}k_k\bar{x}^{y_k}=\sum_{k:~y'_k=\eta}k_k\bar{x}^{y_k}.
%\end{equation*}
%If a mass-action kinetic system admits a complex balanced equilibrium, all of its equilibria are complex balanced. The kinetic system admitting a complex balanced equilibrium is called a complex balanced system. \cite{Horn1972} proved that each non-negative stoichiometric compatibility class contains a unique complex balanced equilibrium for a complex balanced system.
\cite{G2018b} further proved the existence, uniqueness, and local asymptotic stability of equilibrium in DeCBMAS relative to its stoichiometric compatibility class utilizing the following Lyapunov-Krasovskii functional $V:\mathscr{C}_{+}\rightarrow \mathbb{R}_{\geq 0}$, given by  
 \begin{small}
 \begin{equation}\label{eq:Vd}
\begin{split}
&V(\psi)=\sum_{j=1}^{n}(\psi_{j}(0)(\ln(\psi_{j}(0))-\ln(\bar{x}_{j})-1)+\bar{x}_{j})\\
&+\sum_{i=1}^{r}k_{i}\int^{0}_{-\tau_{i}}\left\{(\psi (s))^{y_{\cdot i}}\left[\ln((\psi(s)^{y_{\cdot i}})-\ln(\bar{x}^{y_{\cdot i}})-1\right]+\bar{x}^{y_{\cdot i}}\right\} ds.
\end{split}
\end{equation}
\end{small}

Some definitions and notations for persistence of DeMAS are given in the following part which was firstly proposed by \cite{H2019} and its non-delayed version can be found in \cite{Anderson2008}. 
\begin{definition}\em(persistence)
An $\mathcal{M}_{\mathrm{De}}=(\mathcal{S,C,R}, k, \tau)$ described by (\ref{eq:dde}) is persistence if any forward trajectory $x^{\psi}(t)\in\mathbb{R}^{n}_{\geq 0}$ with a positive initial condition $\psi\in \mathscr{C}_{+}$ satisfies 
\begin{equation}
\liminf_{t\rightarrow \infty}x^{\psi}_{j}(t)>0~\mathrm{for}~\mathrm{all}~j\in\{1,\cdots,n\}.
\end{equation}
\end{definition}
As bounded systems, the persistence of DeCBMASs, can be reflected by $\omega$-limit set. 
\begin{definition}\em($\omega$-limit set)
The $\omega$-$limit$ set for the trajectory $x^{\psi}(t)$ with a positive initial condition $\psi\in \mathscr{C}_{+}$ is
\begin{equation*}
\begin{split}
    \omega(\psi):=&\{\phi\in \bar{\mathscr{C}}_{+}~|~x^{\psi}_{t_{N}}\rightarrow \phi,~\text{for some time sequence} \\&~t_{N}\rightarrow \infty~
    \text{with} ~t_{N}\in\mathbb{R}\}.
    \end{split}
\end{equation*}
The $\omega$-limit set is actually a positive invariant set of the corresponding trajectory.
\end{definition}

\begin{definition}\em(persistence for bounded trajectories)\label{def:FW}
For a DeMAS $\mathcal{M}_{\mathrm{De}}=(\mathcal{S,C,R},k,
\tau)$ with trajectories bounded , it is persistent if
\begin{equation}\label{eq:wL}
 \omega(\psi)\cap(\cup_{W}L_W)=\emptyset
 , ~~~~\forall~\psi\in\mathscr{C}_{+},
 \end{equation} 
where $W$ is a non-empty subset of $\mathcal{S}$ and
\begin{equation}\label{def:Lw}
   L_W=\left\{\theta\in \bar{\mathscr{C}}_{+}\big|\substack{\theta_j(s)=0,~ X_j\in W,\\\theta_j(s)\neq 0,~X_j\notin W,}~~\forall s\in [-\tau_{\rm{m}},0]\right\}
\end{equation}
is called a boundary of  $\bar{\mathscr{C}}_{+}$ and $\cup_{W}L_W$ denotes the set of all boundaries of $\bar{\mathscr{C}}_{+}$. Further $F_W=L_W\cap D_{\psi}$ denotes a boundary of the stoichiometric compatibility class $D_{\psi}$.
\end{definition}
As $\omega(\psi)\in\mathcal{D}_{\psi}$, \eqref{eq:wL} means $\omega(\psi)\cap(\cup_WF_W)=\emptyset$. Actually, some boundaries $F_W$ have already been excluded from the possibility of containing $\omega$-limit points, such as non-semilocking boundary, vertex and facet which are defined as follows.

\begin{definition}\em(semilocking boundary)
For an $\mathcal{N}=(\mathcal{S,C,R})$, a non-empty symbol set $W\subseteq\mathcal{S}$ is called a semilocking set if it satisfies: $W\cap\mathrm{supp}~y_{\cdot i}\neq \emptyset$ if $W\cap\mathrm{supp}~y'_{\cdot i}\neq\emptyset$ for each $R_i:y_{.i}\to y'_{.i}$. The corresponding boundary $F_W$ is called a semilocking boundary. Further, $W$ is called the trivial semilocking set if it is $\mathcal{S}$; otherwise called nontrivial semilocking set. 
 \end{definition}

\begin{definition}[\cite{Zhang2021}] Consider an $\mathcal{M}_{\mathrm{De}}=(\mathcal{S,C,R}, k,\tau)$ described by (\ref{eq:dde}) with an initial function $\psi\in\bar{\mathscr{C}_{+}}$. $F_{W}$ defined in Definition \ref{def:FW} is a facet (vertex) if the dimension of $\pi_{W^c}(\mathscr{S})\triangleq \{\pi_{W^c}(v)\vert v\in \mathscr{S}\}$ is $\dim\mathcal{S}-1$ (zero) where  $W^c\triangleq\mathcal{S}\backslash W$ is the complement set of $W$ and $\pi_{W^c}(\mathscr{S})$ means the set of the projection of each  $v\in\mathscr{S}$ on $W^c$. Accordingly $W$ is called the facet-type (vertex-type) set. 
\end{definition}

%%%%%%%%%%%%%%%%%%%%%%%%%%%%%%%%%%%%%%%%
\section{Problem Motivation}\label{sec:3}
\begin{figure*}[htbp]
    \centering
    \includegraphics[height=4cm,width=18cm]{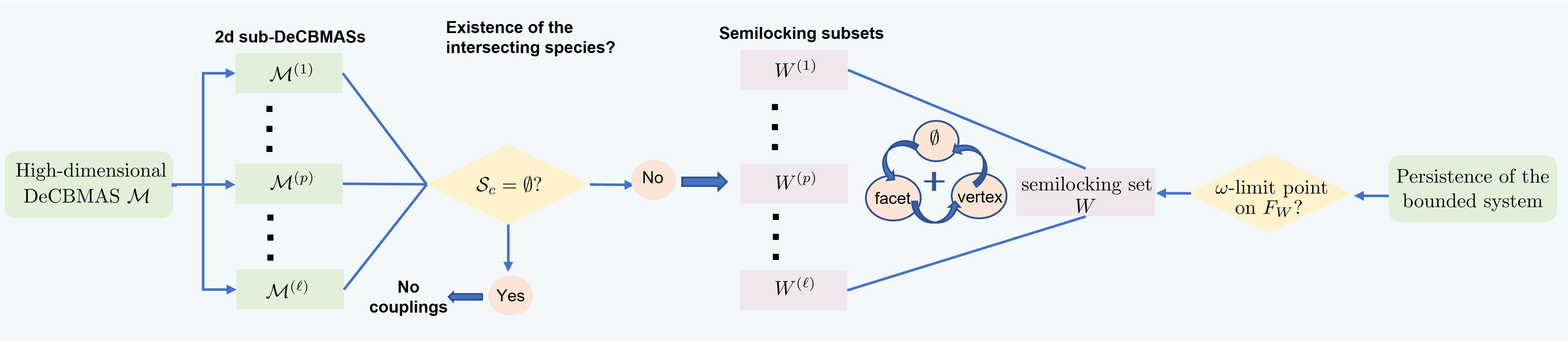}
    \caption{The logical graph of the Section \ref{sec:3}.}
    \label{fig:31}
\end{figure*}
 \cite{Zhang2021} has demonstrated the persistence of $2$-dimensional ($2d$) DeCBMASs ($\dim{\mathscr{S}}=2$) by thoroughly investigating their various semilocking sets (each $W$ can only be a vertex semilocking set or a facet semilocking set). However, conducting a persistence analysis for higher-dimensional DeCBMASs, such as the following example, is challenging due to the intricate conditions governing their boundaries. 
\begin{example}\label{ex:1}
Consider the following $6d$ DeCBMAS
\begin{figure}[H]
    \centering
    \includegraphics[height=5cm,width=8.5cm]{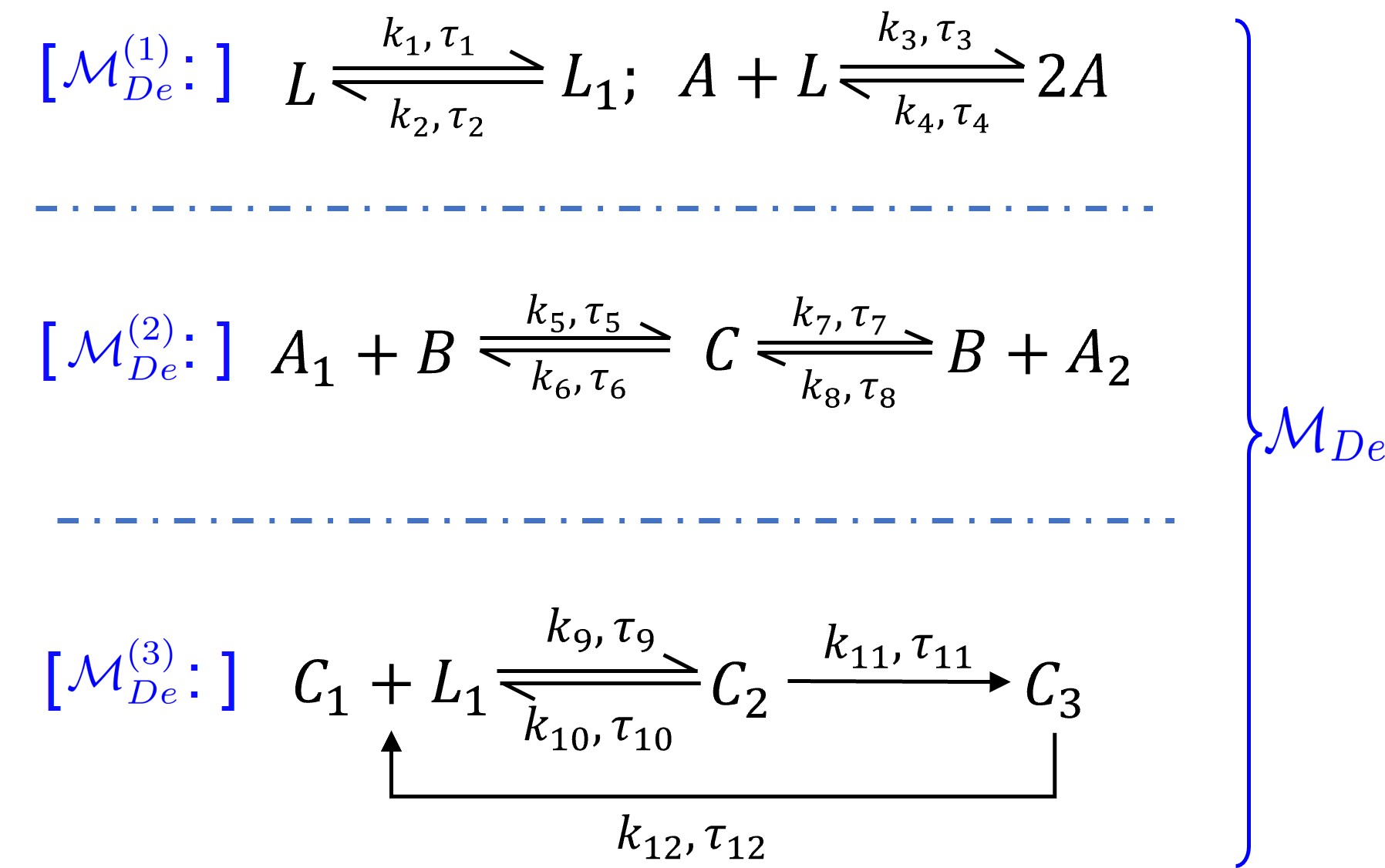}
    \caption{The Combined DeCBMAS $\mathcal{M}_{\mathrm{De}}$.}
    \label{fig:ex1}
\end{figure}
    %\begin{equation*}
%\left. \begin{matrix}
 %	\mathcal{M}_{\mathrm{De}}^{(1)}:& \xymatrix{L\ar @{ -^{>}}^{~\tau_1,~k_1}@< 1pt> [r]& L_1 \ar  @{ -^{>}}^{~\tau_2,~k_2}  @< 1pt> [l]}\\                           
  %&\xymatrix{A+L\ar @{ -^{>}}^{~\tau_3,~k_3}@< 1pt> [r]& 2A \ar  @{ -^{>}}^{~\tau_4,~k_4}  @< 1pt> [l]}\\
   %\mathcal{M}_{\mathrm{De}}^{(2)}:& \xymatrix{A_1+B \ar @{ -^{>}}^{~~\tau_5,~k_5}  @< 1pt> [r] & C \ar @{ -^{>}}^{~~\tau_6,~k_6} @< 1pt> [l] \ar @{ -^{>}}^{\tau_7,k_7~~}  @< 1pt> [r] & B+A_2 \ar @{ -^{>}}^{\tau_8,k_8~~} @< 1pt> [l]}\\
   %\mathcal{M}_{\mathrm{De}}^{(3)}:&\xymatrix{C_1+L_1\ar @{ -^{>}}^{~\tau_9,~k_9}@< 1pt> [r]& C_2 \ar  @{ -^{>}}^{~\tau_{10},~k_{10}}  @< 1pt> [l] \ar @{ -^{>}}^{\tau_{11},k_{11}~~}  @< 1pt> [r] & C_3 \ar `d/15pt[l] `[ll]_{\tau_{12},k_{12}} [ll]}
%\end{matrix}\right \} \mathcal{M}_{\mathrm{De}}
%\end{equation*}

Although the semilocking boundaries of $\mathcal{M}_{De}$ are complex, they can be understood as a combination of three subsets corresponding to three sub-systems $\mathcal{M}_{\mathrm{De}}^{(p)}, p=1,2,3$. Each of these sub-systems $\mathcal{M}_{\mathrm{De}}^{(p)}$ can be described as a $2d$ DeCBMAS. The first pair of reactions represent the interconversion between two states of the biomolecule denoted as $L$ and $L_1$. $A$ and $A_1$ are the same species, the system combined of the second pair of reactions in $\mathcal{M}_{\mathrm{De}}^{(1)}$ and $\mathcal{M}_{De}^{(2)}$ is actually the Edelstein system \citep{E1970} describing the typical enzyme catalysis reaction system in biochemical process, where $A$, $B$, $A_2$, $C$ denote substrate, enzyme, product and intermediate, respectively. And $\mathcal{M}_{\mathrm{De}}^{(3)}$ is actually the T-cell kinetics-proofreading process. 
\end{example}
 A naive idea is to combine $2d$ complex balanced systems to capture the persistence of higher-dimensional DeCBMASs like Example \ref{ex:1}. Towards this purpose, we assume a DeCBMAS of $\mathcal{M}_{\mathrm{De}}=(\mathcal{S,C,R}, k, \tau)$ is a combination of some $2d$ DeCBMASs $\mathcal{M}_{\mathrm{De}}^{(p)}\triangleq(\mathcal{S}^{(p)},\mathcal{C}^{(p)},\mathcal{R}^{(p)}, k^{(p)}, \tau^{(p)}),~p=1,\cdots,\ell$.
For each $\mathcal{M}_{\mathrm{De}}^{(p)}$, suppose it contains $n_p$ species and $r_p$ reactions, and $y_{{.i}^{(p)}}$ and $y'_{{.i}^{(p)}}$ are the reactant and product complexes of the $i$-th reaction in it, respectively. Then, we have 
	\begin{equation*}
	\begin{split}
    &	\mathcal{S}=\bigcup_{p=1}^{\ell}\mathcal{S}^{(p)}; ~\mathcal{C}=\bigcup_{p=1}^{\ell}\bigcup_{i=1}^{r_p} \{y_{.i}^{(p)},{y'_{.i}}^{(p)}\};\\
&\mathcal{R}=\bigcup_{p=1}^{\ell}\mathcal{R}^{(p)}=\bigcup_{p=1}^{\ell}\bigcup_{i=1}^{r_p}\{y_{.i}^{(p)}\to {y'_{.i}}^{(p)}\};\\
	&k=\bigcup_{p=1}^{\ell}k^{(p)}=\bigcup_{p=1}^{\ell}\bigcup_{i=1}^{r_p}\{k_{i}^{(p)}\}; \tau=\bigcup_{p=1}^{\ell}\tau^{(p)}=\bigcup_{p=1}^{\ell}\bigcup_{i=1}^{r_p}\{\tau_i^{(p)}\},
	\end{split}
	\end{equation*}
where 
	\begin{equation*}
	y_{ji}^{(p)}=
	\left\{
     \begin{matrix}
     	y_{{ji}^{(p)}},&\text{if}~S_j\in\mathcal{S}^{(p)}\\
     0,&\text{if}~S_j\notin\mathcal{S}^{(p)}
     \end{matrix}
     \right. ;	{y'_{ji}}^{(p)}=
	\left\{
     \begin{matrix}
     	y'_{{ji}^{(p)}},&\text{if}~S_j\in\mathcal{S}^{(p)}\\
     0,&\text{if}~S_j\notin\mathcal{S}^{(p)}.
     \end{matrix}
     \right.
	\end{equation*}
 As the persistence of $\mathcal{M}_{\mathrm{De}}$ can be reflected by 
 $\omega(\psi)\cap(\cup_{W}F_W)=\emptyset$
for any $\psi\in\mathscr{C}_{+}$ and any semilocking boundary $F_W$, we can divide each semilocking set $W$ of $\mathcal{M}_{\mathrm{De}}$ into $\ell$ subsets $W^{(p)}\triangleq W\cap\mathcal{S}^{(p)}, p=1,\cdots, \ell$. Moreover, we have 
 \begin{lemma}[\cite{Zhang2023}]
     For an $\mathcal{N}=(\mathcal{S,C,R})$ composed of a finite subnetworks $\mathcal{N}^{(p)}=(\mathcal{S}^{(p)},\mathcal{C}^{(p)},\mathcal{R}^{(p)}), p=1,\cdots,\ell$, if $W\subseteq \mathcal{S}$ is a semilocking set of $\mathcal{N}$ and $W^{(p)}\triangleq W\cap\mathcal{S}^{(p)}$ is nonempty, then $W^{(p)}$ is a semilocking set in $\mathcal{N}^{(p)}$. 
 \end{lemma}
 
\begin{remark}
The above lemma indicates that when $\mathcal{M}_{\mathrm{De}}^{(p)}$ is a $2d$ sub-DeCBMASs, $W^{(p)}$ corresponds to an empty set, or a vertex-type, or a facet-type semilocking set associated with $\mathcal{M}_{\mathrm{De}}^{(p)}$. We thus classify $W$ into three cases according to $W^{(p)}$:
\begin{itemize}
    \item {Case I:} $\exists~ p$, $W^{(p)}$ is a facet-type semilocking set;
    \item {Case II:} $\forall~ p$, $W^{(p)}$ is a vertex-type semilocking set; 
    \item Case III: $\forall~p$, $W^{(p)}$ is a vertex-type semilocking set or an empty set.
\end{itemize}
%Each $W$ is accordingly called the FT semilocking set, the VeT semilocking set, and the VeET semilocking set, respectively. 
\end{remark}

We use the following example to make an illustration.  
\begin{example}\label{ex:2}
$\mathcal{M}_{\mathrm{De}}$ is a 5d DeCBMAS in the form of
\begin{figure}[H]
    \centering
    \includegraphics[height=6cm,width=8cm]{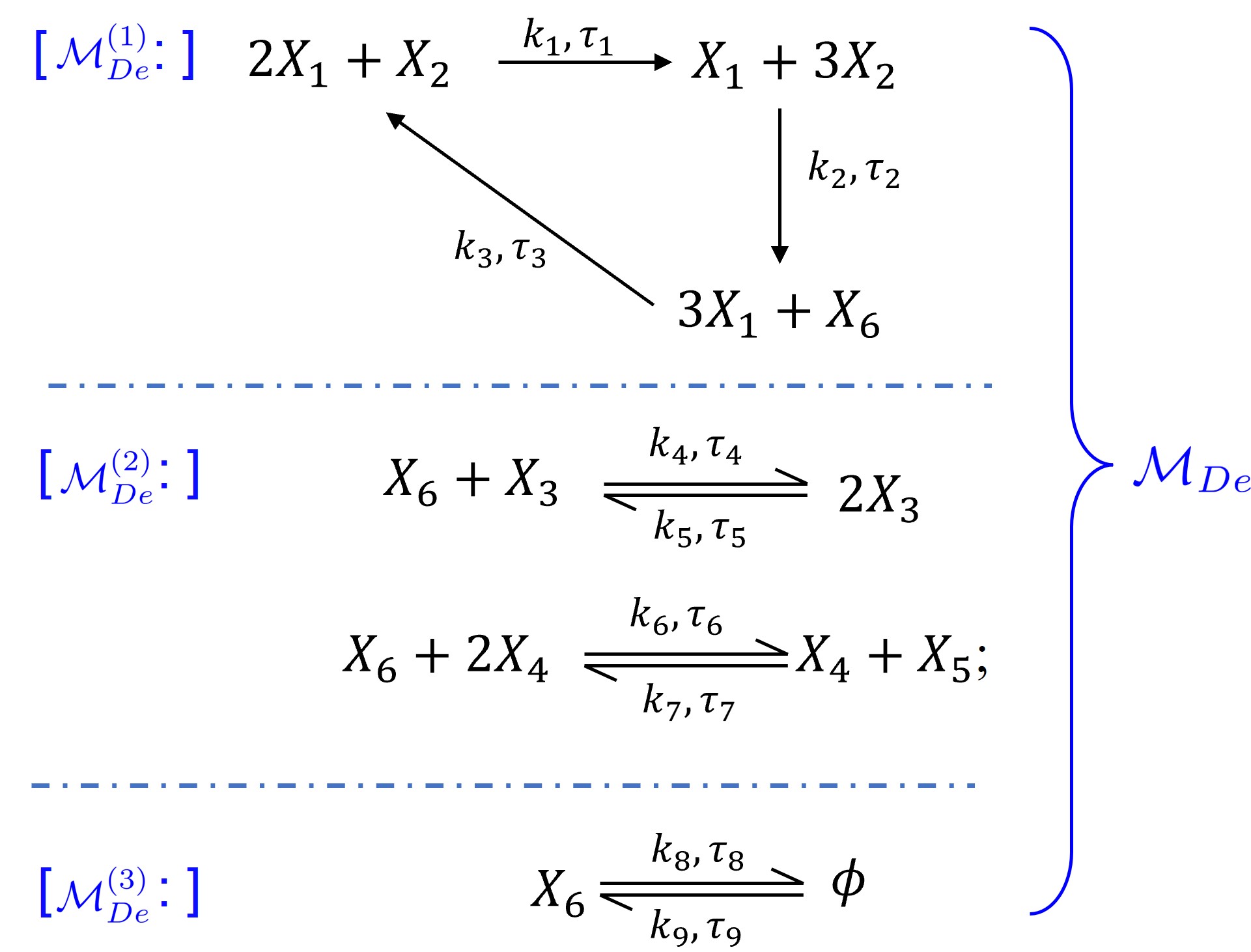}
    \caption{The Combined DeCBMAS $\mathcal{M}_{\mathrm{De}}$.}
    \label{fig:ex2}
\end{figure}
%\begin{equation*}
%\left. 
%\begin{matrix}
%\mathcal{M}_{\mathrm{De}}^{(1)}:&\xymatrix{2X_1+X_2\ar ^{\tau_1,~k_1~~} [r] &X_{1}+3X_{2} \ar ^{~\tau_2,~k_2} [d]\\ 
%&3X_{1}+X_6\ar ^{~\tau_3,~k_3~~} [lu]} \\ 
%\hdashline \mathcal{M}_{\mathrm{De}}^{(2)}:&\xymatrix{X_6+X_3\ar @{ -^{>}}^{~\tau_4,~k_4}@< 1pt> [r]& 2X_3 \ar  @{ -^{>}}^{~\tau_5,~k_5}  @< 1pt> [l]}\\ 
 %&\xymatrix{X_6+2X_4\ar @{ -^{>}}^{~\tau_6,~k_6}@< 1pt> [r]& X_4+X_5 \ar  @{ -^{>}}^{~\tau_7,~k_7}  @< 1pt> [l]}\\
 %\hdashline
 %\mathcal{M}_{\mathrm{De}}^{(3)}:&
  %   \xymatrix{X_6\ar @{ -^{>}}^{~\tau_8,~k_8}@< 1pt> [r]& \emptyset \ar  @{ -^{>}}^{~\tau_9,~k_9}  @< 1pt> [l]}
%\end{matrix}\right \} \begin{matrix}
 %   & \mathcal{M}_{\mathrm{De}}
    %~\mathrm{with}   \\    &\mathcal{S}_c=\{X_\}.
%\end{matrix}
%\end{equation*}

 The semilocking sets of $\mathcal{M}_{\mathrm{De}}$ are listed as follows:
\begin{equation*}
\begin{split}
& W_1=\{X_1\}; W_2=\{X_3\}; W_3=\{X_4\}; W_4=\{X_1, X_2\};\\
&W_5=\{X_1, X_3\}, W_6=\{X_1, X_2, X_3\}; W_7=\{X_1, X_4\};\\& W_8=\{X_1, X_2, X_4\}; W_{9}=\{X_1, X_3, X_4\};\\
& W_{10}=\{X_3, X_4\}; W_{11}=\{X_3, X_4, X_5\}\\
& W_{12}=\{X_1, X_4, X_5\}; W_{13}=\{X_1, X_2, X_3, X_4\};\\
& W_{14}=\{X_1, X_2, X_4, X_5\}; W_{15}=\{X_1, X_3, X_4, X_5\};
\\ & W_{16}=\{X_1, X_2, X_3, X_4, X_5\}.
\end{split}
\end{equation*}
Clearly $W_1$, $W_2$ and $W_3$ belong to Case I. For $W_5$, since $W_5^{(1)}=\{X_1\}$ and $W_5^{(2)}=\{X_3\}$ are both facet-type semilocking sets, it belong to \textit{Case I}, too. The similar analysis may apply to $W_6$, $W_7$, $W_8$, $W_9$, $W_{12}$, $W_{14}$ and $W_{15}$ to indicate them also in \textit{Case I}. A look at $W_4$ might reveal that $W_4^{(1)}=\{X_1, X_2\}$, $W_4^{(2)}=\emptyset$ and $W_4^{(3)}=\emptyset$. As $W_4^{(1)}$ is a vertex-type semilocking set of $\mathcal{M}^{(1)}$, $W_4$ falls into \textit{Case III}. Similarly $W_{10}$, $W_{11}$, $W_{13}$ and $W_{16}$ are all in \textit{Case III}. In this network, there is no semilocking set in \textit{Case II}. However for the closed form of $\mathcal{M}_{\mathrm{De}}$ denoted as $\mathcal{M'}_{\mathrm{De}}$ derived by removing $\mathcal{M}_{\mathrm{De}}^{(3)}$ from $\mathcal{M}_{\mathrm{De}}$, namely, $X_j$ has no input or output from the external environment, $W_{13}$ is also a semilocking set of $\mathcal{M'}_{\mathrm{De}}$ and would fall in \textit{Case II}.

%, which can be solved using the method described in \cite{Zhang2021}. Although the remaining semilocking sets do not fall into the categories of vertex-type or facet-type sets, they can still be classified into the three cases mentioned earlier by considering $\mathcal{M}_{\mathrm{De}}$ as a combination of three subsystems, denoted as $\mathcal{M}^{(p)}_{\mathrm{De}}$ where $p=1, 2, 3$.
\end{example}

In this work, it is expected the combined higher dimensional $\mathcal{M}_{\mathrm{De}}$ can inherit the persistence of $2d$ sub-DeCBMASs.
Apart from the boundary characteristics of $\mathcal{M}_{\mathrm{De}}$ influenced by individual $2d$ subsystems, the manner in which they are networked together is also crucial in determining whether $\mathcal{M}_{\mathrm{De}}$ can inherit the persistence of the subsystems. The analysis method, designed for this purpose, has to be able to handle the coupling of dynamics induced by intersecting species between sub-DeCBMASs. We use $\mathcal{S}_c$ to represent the intersecting species set in which every species appears in more than one subsystems. Clearly the trivial scenario occurs when $\mathcal{S}_c=\emptyset$, meaning that $\mathcal{S}^{(p_1)}\cap \mathcal{S}^{(p_2)}=\emptyset$ for any two $p_1$ and $p_2$. In this case, the persistence of $\mathcal{M}_{\mathrm{De}}$ is evident, as each $2d$ $\mathcal{M}_{\mathrm{De}}^{(p)}$ is independently persistent and no couplings between subsystems. Example \ref{ex:1} illustrates that when $A$, $A_1$ and $C$, $C_1$, $C_2$ represent distinct species, the persistence of the higher-dimensional system $\mathcal{M}_{\mathrm{De}}$ can be derived directly from the subsystems.
Hence, the primary focus of this paper is on the persistence of higher-dimensional DeCBMASs composed of $2d$ DeCBMASs in the case of $\mathcal{S}_c\neq \emptyset$. In this scenario, the persistence of $\mathcal{M}_{\mathrm{De}}$ is difficult to capture due to the presence of dynamics coupling between subsystems. For instance, there may be instances where $A$, $A_1$ or $C$, $C_1$, $C_2$ correspond to the same species in Example \ref{ex:1}.

To tackle this issue, a possible solution is to figure out the relation between nonempty $\mathcal{S}_c$ and the semilocking set of the combined DeCBMAS, since the latter plays an important role in diagnosing persistence. We thus discuss two possibilities according to their relations: one is that there is no intersecting species in the nontrivial semilocking set, and the other is that there is intersecting species in the nontrivial semilocking set. Figure \ref{fig:31} exhibits the logical structure of our whole idea.

\section{Higher-dimensional DeCBMAS composed of $2d$ sub-DeCBMASs with intersecting species in non-semilocking set}\label{sec:4}
In this section, we consider the first possibility that there is no intersecting species in the nontrivial semilocking set, i.e, for any nontrivial semilocking set $W$ in the combined DeCBMAS of $\mathcal{M}_{\mathrm{De}}=(\mathcal{S,C,R},k,
\tau)$, $\mathcal{S}_c\cap W=\emptyset$. The dynamics coupling only occurs to the species in non-semilocking sets of $\mathcal{M}_{\mathrm{De}}$. Like Example \ref{ex:2}, $\mathcal{S}_{c}=\{X_6\}$ where the intersecting species $X_6$ only exists in non-semilocking sets of $\mathcal{M}_{\mathrm{De}}$. The presence of intersecting species necessitates the development of a new approach for persistence analysis. Based on the cases of semilocking sets of the combined DeCBMAS, the corresponding boundaries $F_W$ must be also considered from the above three cases. And the persistence of $\mathcal{M}_{\mathrm{De}}$ will be conducted by dealing with each of the three possible $F_W$. Note that we will put all proofs of results in this section and others to Appendix to improve readability.

\subsection{Semilocking set in Case I}
The following part is devoted to the studying of the case I, namely, some $W^{(p)}$ is a facet-type semilocking set. Typical examples of this kind of semilocking set are $W_5, W_6$.
\begin{lemma}[\cite{Zhang2023cdc}]\label{lem:c3}
  A DeCBMAs denoted as $\mathcal{M}_{\mathrm{De}}=(\mathcal{S,C,R},k,\tau)$ is composed of a finite number of DeCBMASs $\mathcal{M}_{\mathrm{De}}^{(p)}, p=1,\cdots,\ell$. $F_W$ does not exist $\omega$-limit points of each trajectory starting from $\psi\in\mathscr{C}_+$, namely, $\omega(\psi)\cap F_W=\emptyset$, if  $\mathcal{S}_c\cap W=\emptyset$ and $W^{(p)}$ is a facet-type semilocking set of $\mathcal{M}^{(p)}$.
\end{lemma}
\begin{remark}
Based on the insights gained from the proof of Lemma \ref{lem:c3}, it can be observed that when $S_c\cap W=\emptyset$, the concentration of each species $X_j\in W^{(p)}$ in $\mathcal{M}_{De}$ cannot reach 0  independent of the forms of other subnetworks as long as $W^{(p)}$ is a facet-type semilocking set of $\mathcal{M}^{(p)}$. 
\end{remark}

Then we can generalize Lemma \ref{lem:c3} into the following Lemma.
\begin{lemma}[Boundary of Case I]\label{thm:4.9}
A DeCBMAS $\mathcal{M}_{\mathrm{De}}=(\mathcal{S, C, R},k,\tau)$ composed of $\ell$ sub-DeCBMASs $\mathcal{M}^{(p)}, p=1,\cdots,\ell$ is with dynamics of the form (\ref{eq:dde}). Let $W$ be a semilocking set, then $F_W\cap \omega(\psi)=\emptyset$ for each $\psi\in\mathscr{C}_+$ if its corresponding semilocking set $W$ is in Case I and satisfies $\mathcal{S}_c\cap W=\emptyset$.
\end{lemma}
Naturally, the persistence of some high-dimensional DeCBMASs with only Case I semilocking sets are derived.
\begin{theorem}
 $\mathcal{M}_{\mathrm{De}}=(\mathcal{S,C,R},k,\tau)$ is a DeCBMAS which can be seen as a combination of a finite number of sub-DeCBMASs $\mathcal{M}^{(p)}_{\mathrm{De}}, p=1,\cdots, \ell$. $\mathcal{M}_{\mathrm{De}}$ is persistent if each its semilocking set $W$ satisfies $W$ in Case I and $W\cap\mathcal{S}_c=\emptyset$.
\end{theorem}

%Also it holds for some DeMASs that are not complex balanced.
%\begin{corollary}\label{lem:c3c}
 %   $\mathcal{M}_{\mathrm{De}}$ is a conservative DeMAS which can be seen as the combination of $\ell$ mass action systems  $\mathcal{M}_{\mathrm{De}}^{(p)}, p=1,\cdots,\ell$. For each semilocking set $W$, if there exists some $p$ such that $\mathcal{M}_{\mathrm{De}}^{(p)}$ is a DeCBMAS and $W^{(p)}$ is the facet-type semilocking set of $\mathcal{M}_{\mathrm{De}}^{(p)}$, then $\mathcal{M}_{\mathrm{De}}$ is persistent. 
%\end{corollary}
%\begin{example}
%$\mathcal{M}_{\mathrm{De}}$ is a typical delayed Com-$\ell$ts-Autoca system defined in \cite{Lu2022}
%\begin{figure}[H]
 %   \centering
   % \includegraphics[height=6cm,width=8cm]{figex3.jpg}
    %\caption{The Combined DeCBMAS $\mathcal{M}_{\mathrm{De}}$.}
    %\label{fig:ex3}
%\end{figure}
%\begin{equation*}
%\left.
%\begin{matrix}
%\begin{array}{cc}
%\begin{matrix}
 %\mathcal{M}_{\mathrm{De}}^{(1)}:&\xymatrix{X_1\ar @{ -^{>}}^{\tau_1,~k_1}@< 1pt> [r]& X_2\ar  @{ -^{>}}^{\tau_2,~k_2}  @< 1pt> [l]}\\
%&\xymatrix{2X_1\ar ^{\tau_3,~k_3~~} [r] &X_{1}+X_{2} \ar ^{~\tau_4,~k_4} [d]\\ 
%& 2X_{2}+X_4\ar ^{~\tau_5,~k_5~~} [lu]}\\
%\end{matrix}\\ 
%\hdashline
%\begin{matrix}
%\mathcal{M}_{\mathrm{De}}^{(2)}:
 %\xymatrix{X_4\ar @{ -^{>}}^{\tau_6,~k_6}@< 1pt> [r]& X_3\ar  @{ -^{>}}^{\tau_7,~k_7}  @< 1pt> [l]}\\
  %   X_4+(m_i-1)X_3\xrightarrow{\tau_{m_i+6}, k_{m_i+6}} m_iX_3\\
   %  m_i=2,3,4
%\end{matrix}
%\end{array}\\
%\end{matrix}
%\right \}
%\begin{matrix}
 %   & \mathcal{M}_{\mathrm{De}}
  %  ~\mathrm{with}   \\    &\mathcal{S}_c=\{X_4\}.
%\end{matrix}
%\end{equation*}

\subsection{Semilocking set in Case II}
This subsection primarily focuses on the investigation of  $F_W$ of $\mathcal{M}_{\mathrm{De}}$ where the corresponding semilocking set $W$ falls into Case II.                                
\begin{lemma}[Boundary of Case II]\label{lem:vv} 
 Let $\mathcal{M}_{\mathrm{De}}=(\mathcal{S,C,R},k,\tau)$ be a DeCBMAS composed of finite number of sub-DeCBMASs $\mathcal{M}_{\mathrm{De}}^{(p)}=(\mathcal{S}^{(p)},\mathcal{C}^{(p)},\mathcal{R}^{(p)},k,\tau),~p=1,\cdots,\ell$. Suppose $W$ is a semilocking set of $\mathcal{M}_{\mathrm{De}}$ in Case II satisfying $\mathcal{S}_c\cap W=\emptyset$. Then $W$ is also a vertex-type semilocking set of $\mathcal{M}_{\mathrm{De}}$ and there does not exist an $\omega$-limit point of the trajectory starting from any $\psi\in \mathscr{C}_+$ on $F_W$.
\end{lemma}
The following example indeed provides supporting evidence for Lemma \ref{lem:vv} and its proof. 
\begin{example}\label{ex:4.2}
    Consider the closed form system $\mathcal{M'}_{\mathrm{De}}$ which is composed of the first two subsystems $\mathcal{M}^{(1)}$ and $\mathcal{M}^{(2)}$ in Example \ref{ex:2}. $W_{13}$ is the semilocking set in the Case II of $\mathcal{M'}_{\mathrm{De}}$.  Its corresponding boundary $F_{W_{13}}$ given by $F_{W_{13}}=L_{W_{13}}\cap\mathcal{D}_{\psi}$ where 
    $$L_{W_{13}}=\{\theta\in\bar{\mathscr{C}}_+\vert (0,0,0,0,\theta_5, \theta_6), \theta_5(s)>0, \theta_6(s)>0\}$$
    is still a vertex of each $\mathcal{D}_{\psi}$ of $\mathcal{M'}_{\mathrm{De}}$ as there does not exist a vector $v\in \mathscr{S}$ in the form $v=(0,0,0,0,v_5,v_6), v_5\neq 0$ or $v_6\neq 0$. Therefore, the possibility of $F_W$ containing $\omega$-limit points is excluded.
\end{example}
 Then the following results can be concluded directly from the Lemma \ref{lem:vv}.
 \begin{theorem}\label{thm:c1}
     Let $\mathcal{M}_{\mathrm{De}}$ be a DeCBMAS composed of $\ell$ sub-DeCBMASs $\mathcal{M}_{\mathrm{De}}^{(p)}, p=1,\cdots, \ell$. If each semilocking set $W$ is in Case I or Case II and satisfies $W \cap \mathcal{S}_c = \emptyset$, then $\mathcal{M}_{\mathrm{De}}$ is persistent.
 \end{theorem}

 \subsection{Semilocking set in Case III}
This subsection discusses the boundary of the remaining Case III. For each such $W$, without loss of generality, we assume that its subsets can be represented as some vertex-type semilocking sets $W^{(p)}, p=1,\cdots,\ell_1$ of $\mathcal{M}^{(p)}=(\mathcal{S}^{(p)}, \mathcal{C}^{(p)},\mathcal{R}^{(p)},k^{(p)},\tau^{(p)})$  and some empty sets $W^{(p')}, p'=\ell_1+1, \cdots, \ell$ of $\mathcal{M}^{(p')}=(\mathcal{S}^{(p')}, \mathcal{C}^{(p')},\mathcal{R}^{(p')},k^{(p')}, \tau^{(p')})$. 
%The existing of empty set $W^{(p')}$ compared to Case II resulting in $W$ not being a vertex of the $\mathcal{M}_{\mathrm{De}}$. 
The following example illustrates the impact of existing empty subsets compared to the Case II:
\begin{example}\label{ex:4}
    Let's consider the open DeCBMAS $\mathcal{M}_{\mathrm{De}}$ and its closed form $\mathcal{M'}_{\mathrm{De}}$ once again. In this example, $W_{13}$ represents the semilocking set, and $W_{13}^c=\{X_5, X_6\}$ represents the complement of $W_{13}$.
    \begin{itemize}
    \item For the closed form $\mathcal{M'}_{\mathrm{De}}$, none of the reactions involving the species in $W_{13}^c$ can occur on the boundary $F_{W_{13}}$. Furthermore, there does not exist any $v\in\mathscr{S}'$ as described in Example \ref{ex:4.2}. This means that each species $X_j\in W_{13}^c$ is restricted by $W_{13}$.
 
    \item In contrast, for the open system $\mathcal{M}_{\mathrm{De}}$, the empty subset $W_{13}^{(3)}$ introduces reactions in $\mathcal{M}_{De}^{(3)}$ that can take place on $F_{W_{13}}$. As a result, each species in $\mathcal{M}_{\mathrm{De}}^{(3)}$ is not restricted by $W_{13}$. For example, in $\mathcal{M}_{\mathrm{De}}$, $X_6\in W^c_{13}$ is not restricted by $W_{13}$. However, $X_5$ is still restricted like in $\mathcal{M'}_{\mathrm{De}}$ since it is not in $\mathcal{M}_{De}^{(3)}$.
    Let's consider another example with the semilocking set $W_4$ from Example \ref{ex:2}. In this case, both subsets $W_4^{(2)}$ and $W_4^{(3)}$ are empty. Each species in the complement set $W_4^c=\{X_3, X_4,X_5,X_6\}$ is involved in reactions in $\mathcal{M}^{(2)}$ and $\mathcal{M}^{(3)}$. Therefore, all species in $W_4^c$ are not restricted.
    %completely removes the restrictions on $X_6$. There exists non-empty vectors $v=(0,0,0,0,0,v_6)\in \mathscr{S}$. Additionally, the reactions in $\mathcal{M}^{(3)}_{\mathrm{De}}$  can take place on the boundary $F_{W_{13}}$. Note that the concentration of $X_5$ denoted as $x_5$ is still be restricted by the species in $W$ and $x_5$ is a constant on the boundary $F_{W_{13}}$ as the corresponding subsystem $\mathcal{M}^{(3)}_{\mathrm{De}}$ of the empty subset $W^{(3)}_{13}$ does contain the $X_5$. 
    \end{itemize}
\end{example}
 \begin{remark}\label{rek:1}
      Based on the example discussed above, we can conclude that $W^c$ can be divided into two parts according to whether the species is in some $\mathcal{M}^{(p')}$ with empty $W^{(p')}$. Let $W^{cn}\subset W^c$ denote the set of species appearing in $\mathcal{M}_{\mathrm{De}}^{(p')}$ and $W^{cv}\subset W^c$ denote the set of the species only existing in $\mathcal{M}_{\mathrm{De}}^{(p)}$. Specially, each $W$ in Case II satisfies $W^c=W^{cv}$ and $W$ in Case III satisfies $W^{cn}\neq\emptyset$.   
\end{remark}
Thus the analysis of $F_W$ with $W$ in Case III can be conducted based on the potential emptiness of the set $W^{cv}$. 
By utilizing the fact that the concentration of each species in $W^{cv}$ is a constant on $F_W$, we can obtain the following lemma.
\begin{lemma}[Boundary of Case III with $W^{cv}\neq\emptyset$]\label{lem:c21}
Given $\mathcal{M}_{\mathrm{De}}=(\mathcal{S,C,R},k,\tau)$ as a DeCBMAS composed of finite number of sub-DeCBMASs $\mathcal{M}_{\mathrm{De}}^{(p)}, p=1,\cdots, \ell$. There does not exist any $\omega$-limit point of the trajectory with a positive initial point on the boundary $F_W$ if the corresponding $W$ is in Case III with $W^{cv}\neq \emptyset$ and $W\cap\mathcal{S}_c=\emptyset$.
\end{lemma}
In the case where $W^{cv}$ is empty, the similar result also holds by transforming the coupled subsystems of  $\mathcal{M}_{\mathrm{De}}$ into independent generalized delayed mass action systems (GDeMAS) and further analyzing the generalized mass action system.
\begin{figure*}
    \centering
    \includegraphics[height=5cm,width=16cm]{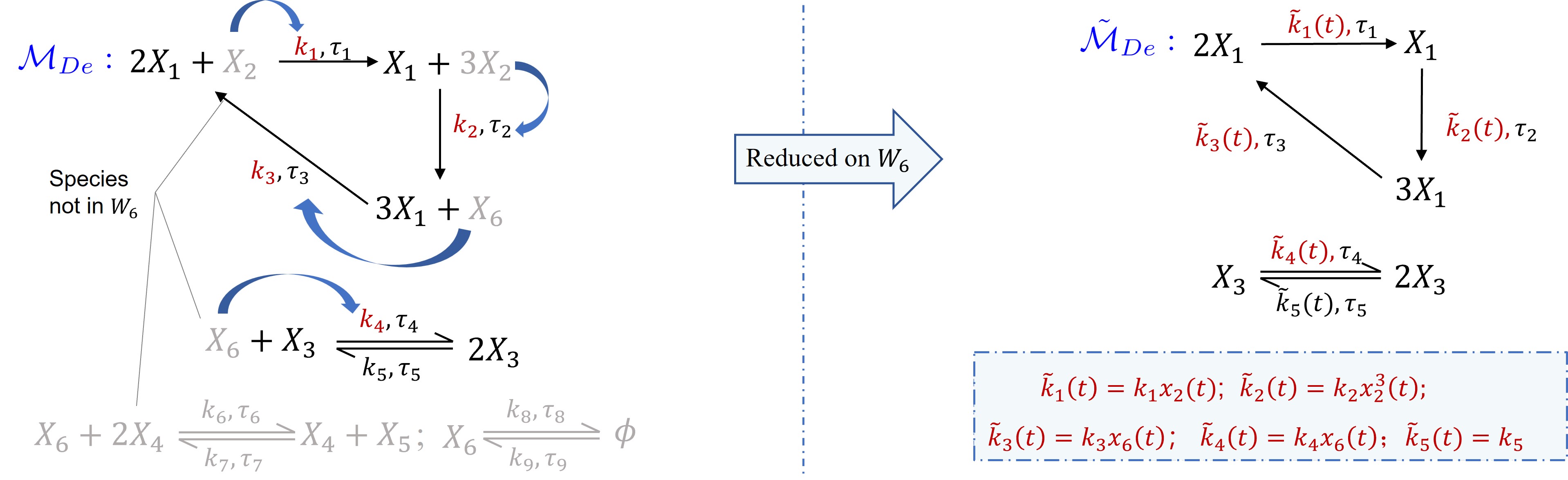}
    \caption{The reduced process of the system $\mathcal{M}_{\mathrm{De}}$ in Example \ref{ex:2} based on $W_6=\{X_1,X_3\}$. The species not in $W_6$ are represented in gray. By reducing, we include them in the kinetic constants $k_i$, thereby transforming $k_i$ into a time-dependent function $k_i(t)$.  Also the reduced system can be derived in the form shown on the right side.}
    \label{fig:3s}
\end{figure*}
\begin{lemma}[Boundary of Case III with $W^{cv}=\emptyset$]\label{lem:c22}
 $\mathcal{M}_{\mathrm{De}}=(\mathcal{S,C,R},k,\tau)$ is a DeCBMAS composed of $\ell$ 2d sub-DeCBMASs $\mathcal{M}_{\mathrm{De}}^{(p)}s$. Then each boundary point on the $F_W$ can not be an $\omega$-limit point of some trajectory starting from the positive point if $W$ is in Case III with $W^{cv}=\emptyset$ and $W\cap\mathcal{S}_c=\emptyset$.
 \end{lemma}
 The reduced system in the delayed case is introduced in the proof of Lemma \ref{lem:c22} to address couplings between subsystems.
\begin{definition}\em(reduced system)\label{def:re}
$\mathcal{M}_{\mathrm{De}}=(\mathcal{S, C, R},k,\tau)$ is a delayed mass action system with dynamics $\dot{x}$. The system $\tilde{\mathcal{M}}_{\mathrm{De}}=(\mathcal{\Tilde{S},\Tilde{C},\Tilde{R}},\tilde{k}(t),\tilde{\tau})$ is called a reduced system of $\mathcal{M}_{\mathrm{De}}$ if
\begin{itemize}      
    \item $\tilde{\mathcal{S}}\subset \mathcal{S}$;  
    \item $\tilde{\mathcal{C}}=\pi_{\mathcal{\tilde{S}}}(\mathcal{C})$, where $\pi_{\mathcal{\tilde{S}}}$ denotes the projection of each complex $y\in \mathcal{C}$ onto components in $\mathcal{\tilde{S}}$. Specially a complex $y\in \mathcal{C}$ with $y_j=0$ for all $X_j\in \tilde{\mathcal{S}}$ corresponds to $\emptyset$ in $\tilde{\mathcal{C}}$, namely $\pi_\mathcal{\tilde{S}}(y)=\emptyset$.
    \item $\mathcal{\tilde{R}}\triangleq \{\pi_{\mathcal{\tilde{S}}}(y_{.i})\rightarrow\pi_{\mathcal{\tilde{S}}}(y'_{.i})\vert  \pi_{\mathcal{\tilde{S}}}(y_{.i})\neq \emptyset~\mathrm{or}~\pi_{\mathcal{\tilde{S}}}(y'_{.i})\neq \emptyset; i=1,\cdots, r\}$.  And each reaction rate parameter of each reaction $\tilde{R}_i: \pi_{\mathcal{\tilde{S}}}(y_{.i})\xrightarrow{\tilde{k}_i(t)}\pi_{\mathcal{\tilde{S}}}(y'_{.i})$ is described as $\tilde{k}_i(t)=k_i\prod_{X_j\notin \tilde{\mathcal{S}}}x_j^{y_{ji}},~\tilde{\tau}_i=\tau_i$. 
    \end{itemize}
\end{definition}
According to the above definition and setting $\Tilde{\mathcal{S}}=W_6$, the reducing process and the corresponding reduced system $\mathcal{\tilde{M}}_{\mathrm{De}}$ of $\mathcal{M}_{\mathrm{De}}$ in Example \ref{ex:2} is shown in the Figure \ref{fig:3s}. We can see that $\mathcal{\tilde{M}}_{\mathrm{De}}$ remains the dynamics of species in $W_6$, but does not consider the dynamics of the intersecting species $X_j$ denoted as $\dot{x}_j, X_j \notin \tilde{\mathcal{S}}$.
Thus through reducing, the system $\tilde{\mathcal{M}}_{\mathrm{De}}$ can be divided into two independent parts without couplings. 

\begin{figure*}
\centering
\subfigure[]{\includegraphics[height=5.8cm,width=8cm]{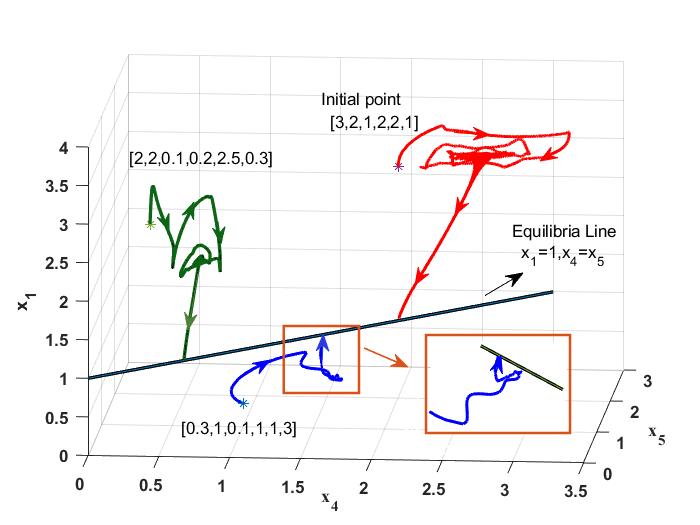}}
\hspace{10mm}
\subfigure[]{\includegraphics[height=5.5cm,width=8cm]{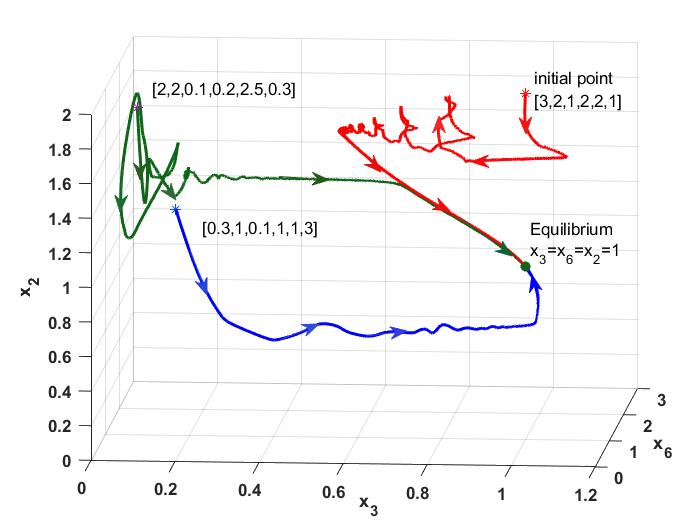}}
\caption{Each equilibrium point $\bar{x}$ of Example \ref{ex:2} is given by $\bar{x}_1=\bar{x}_2=\bar{x}_3=\bar{x}_j=1$ and $\bar{x}_4=\bar{x}_i$. (a) and (b) represent the projections of trajectories in Example \ref{ex:2} with three different initial points and vectors of time delays: $\theta_1, \tau_1$ (in green), $\theta_2, \tau_2$ (in red), and $\theta_3, \tau_3$ (in blue). The values for $\theta_1(s)$, $\theta_2(s)$, and $\theta_3(s)$ are as follows: $\theta_1(s)=[2,2,0.1,0.2,2.5,0.3]$, $\theta_2(s)=[3,2,1,2,2,1]$, and $\theta_3(s)=[0.3,1,0.1,1,1,3]$, where $s\in [-\tau_m, 0]$. The corresponding time delay vectors are $\tau_1=[3,0.3,2,1,3,0.2,5]$, $\tau_2=[1,1,1,2,3,2,0.1]$, and $\tau_3=[2,1,10,1,1,1,0.1]$. In (a), the trajectories are projected into the space $\mathbb{R}^{\mathcal{S}^1}$, where $\mathcal{S}^1=\{X_4,X_i,X_1\}$. In (b), the trajectories are projected into the space $\mathbb{R}^{\mathcal{S}^2}$, where $\mathcal{S}^2=\{S_3, S_j, S_2\}$. }
\label{fig:2.1}
\end{figure*}

\begin{remark}
The dynamical equation $\dot{\tilde{x}}$ of the reduced system $\tilde{\mathcal{M}}_{\mathrm{De}}$ can be expressed as $\dot{\tilde{x}}=\pi_{\tilde{\mathcal{S}}}(\dot{x})$, where $\dot{x}$ is the dynamics corresponding system $\mathcal{M}_{\mathrm{De}}$. 
 Additionally, the reduced system is considered a generalized mass action system. As a result, it maintains the same dynamics for the species within the semilocking set, but transfers the challenges associated with the coupling between subsystems to the complexity of the dynamics within each individual subsystem.
\end{remark}
 %Before this can be achieved, we need the following result.
% \begin{lemma}
 %For a delayed bounded complex  balanced generalized mass action system $\mathcal{\tilde{M}}=(\mathcal{\tilde{S}},\mathcal{\tilde{C}},\mathcal{\tilde{R}}, \tilde{K}, \tilde{\tau})$ where $K=\{k_i(t)\vert i=1,\cdots, r\}$ is the set of reaction rate parameters and each  $M_d^i<k_i(t)<M_u^i$ is a bounded  function with $M_d^i, M_u^i$ all positive constants. $W$ is a semilocking set such that $F_W$ is a facet of the stoichiometric compatibility class, then there does not exist an $\omega$-limit point in $F_W$.
 By considering all possible boundary cases presented in Lemmas \ref{thm:4.9}, \ref{lem:vv}, \ref{lem:c21}, \ref{lem:c22}, the persistence of following $\mathcal{M}_{\mathrm{De}}$ can be obtained.
\begin{theorem}
   A DeCBMAS $\mathcal{M}_{\mathrm{De}}=(\mathcal{S,C,R},k,\tau)$ is composed of finite number of 2d DeCBMASs $\mathcal{M}_{\mathrm{De}}^{(p)}, p=1,\cdots,\ell$. $\mathcal{M}_{\mathrm{De}}$ is persistent if for each non-trivial semilocking set $W$ of $\mathcal{M}_{\mathrm{De}}$, $W\cap \mathcal{S}_c=\emptyset$.
\end{theorem}
 From the above theorem, we can derive that the system $\mathcal{M}_{\mathrm{De}}$ in Example \ref{ex:2} is persistent. Figure \ref{fig:2.1} presents that three trajectories of the system $\mathcal{M}_{\mathrm{De}}$ in Example \ref{ex:2} starting form different initial points and time delays all converge to positive equilibria $x_1=x_2=x_3=x_6=1$, $x_4=x_5$ in the case of $k_i=1, i=1,\cdots, 9$. 
%\section{Some Examples}

%% There are a number of predefined theorem-like environments in
%% ifacconf.cls:
%%
%% \begin{thm} ... \end{thm}            % Theorem
%% \begin{lem} ... \end{lem}            % lem
%% \begin{claim} ... \end{claim}        % Claim
%% \begin{conj} ... \end{conj}          % Conjecture
%% \begin{cor} ... \end{cor}            % Corollary
%% \begin{fact} ... \end{fact}          % Fact
%% \begin{hypo} ... \end{hypo}          % Hypothesis
%% \begin{prop} ... \end{prop}          % Proposition
%% \begin{crit} ... \end{crit}          % Criterion

\section{Higher-dimensional DeCBMAS composed of $2d$ sub-DeCBMASs with the intersecting species in some semilocking sets.}\label{sec:5}
This subsection further explores the scenario where intersecting species also appear in some semilocking set $W$, denoted as $\mathcal{S}_c\cap W\neq \emptyset$ (All proofs of results in this Section can be found in the Section Appendix). 
\begin{lemma}[\cite{Zhang2023cdc}]\label{lem:c531}
$\mathcal{M}_{\mathrm{De}}$ is a DeCBMAS which is a combination of sub-DeCBMASs $\mathcal{M}_{\mathrm{De}}^{(p)}, p=1,\cdots,\ell$ and $W$ is a semilocking set of $\mathcal{M}_{\mathrm{De}}$. $F_W$ cannot contain $\omega$-limit points of some trajectory starting from any $\psi\in\mathscr{C}_{+}$ if each $W^{(p)}$ is a facet-type semilocking set of the subsystem $\mathcal{M}_{\mathrm{De}}^{(p)}$ or an empty set. 
\end{lemma}
With appropriate extensions allowing partial subsets to be vertex-type, the following theorem can be obtained. 
\begin{theorem}[Extension result of Case I]\label{cor:c53}
Suppose $\mathcal{M}_{\mathrm{De}}$ is a DeCBMAS consisting of finite number of sub-DeCBMASs $\mathcal{M}_{\mathrm{De}}^{(p)}$ for $p=1, \cdots, \ell$. $\mathcal{M}_{\mathrm{De}}$ is persistent if for each semilocking set $W\cap \mathcal{S}_c\neq\emptyset$, there exists some species $X_j\in W$ such that $X_j$ is only in facet-type subsets $W^{(p)}$.
\end{theorem}

Thus the persistence of following example can be addressed.
\begin{example}
 Consider a DeCBMAS $\mathcal{M}_{\mathrm{De}}$ in the following form:
\begin{figure}[H]
    \centering
    \includegraphics[height=4.5cm,width=8.5cm]{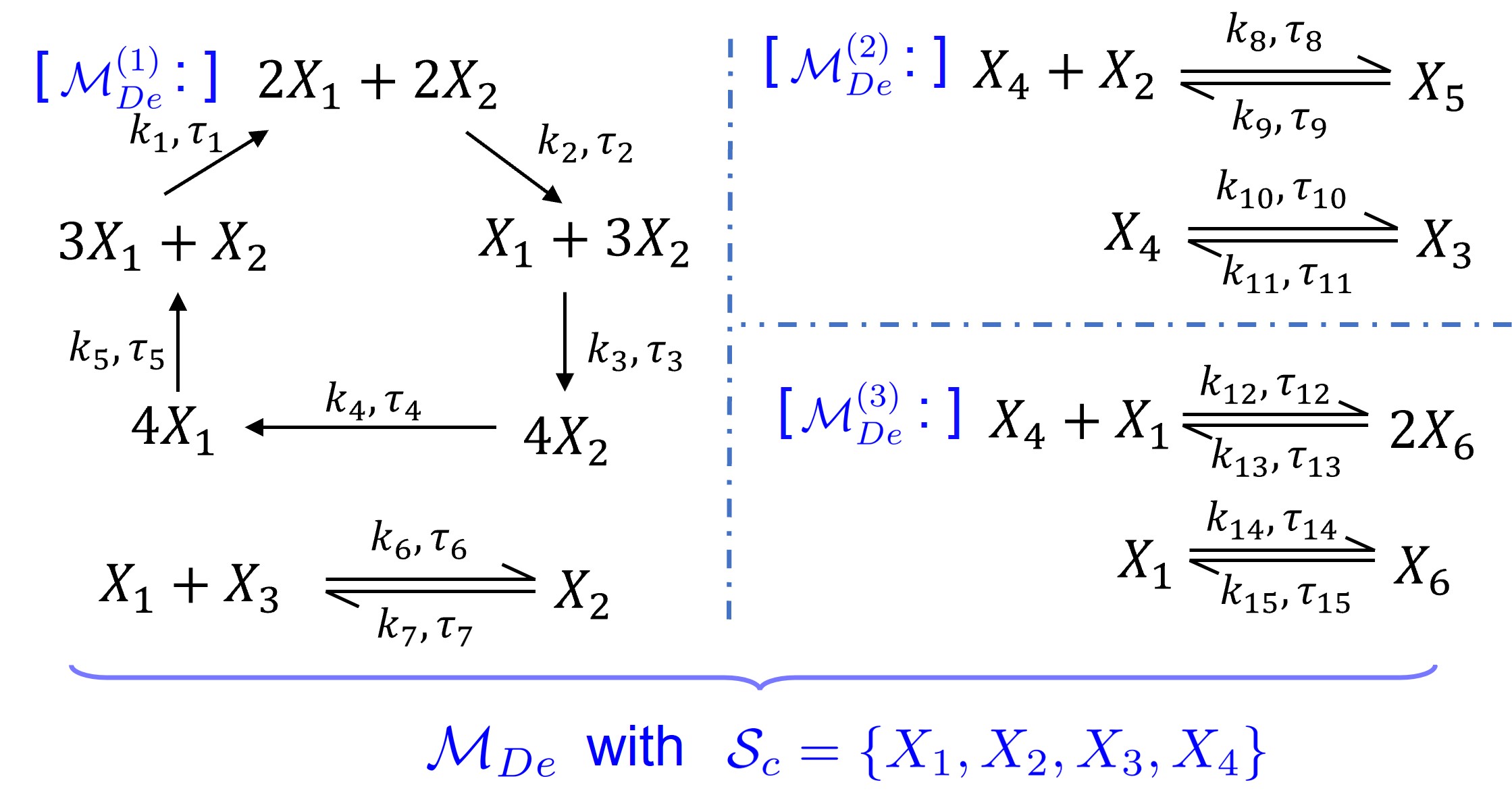}
    \caption{The Combined DeCBMAS $\mathcal{M}_{\mathrm{De}}$.}
    \label{fig:enter-label}
\end{figure}

 $W=\{X_1,X_2,X_5,X_6\}$ is its unique non-trivial semilocking set and corresponding subsets are $W^{(1)}=\{X_1,X_2\}$, $W^{(2)}=\{X_2,X_5\}$, $W^{(3)}=\{X_1,X_6\}$. Thus for this semilocking set, $W\cap\mathcal{S}_c=\{X_1, X_2\}$ and $W^{(1)}$, $W^{(2)}$ are both facet-type semilocking sets while $W^{(3)}$ is a vertex-type semilocking set of their subsystems. Besides, there exist species $X_2, X_5\in W$ that are only contained in the facet-type semilocking subsets. Thus the persistence of high-dimensional DeCBMAS $\mathcal{M}_{\mathrm{De}}$ can be concluded from the Theorem \ref{cor:c53}.
\end{example}

When the condition in Theorem \ref{cor:c53} is not satisfied, it means that each $X_j\in W$ is definitely contained within a vertex-type semilocking set $W^{(p)}$. Similar to the previous chapter, we will conduct an investigation into the elements in $W^c$. In this case, $W^c$ can be divided into three subsets: $W^{tr}$, $W^{sr}$, and $W^{tf}$, defined as follows:
\begin{itemize} 
\item $W^{tf}$: This subset includes species $X_{\bar{j}}$ that are completely unrestricted, meaning they are involved in some reactions that solely consist of species in $W^c$. Thus the concentration of each species in $W^{tf}$ can change even on the boundary $F_W$ due to the occurrence of the reactions.
    \item $W^{sr}$: This subset comprises species $X_{\tilde{j}}$ that are semi-restricted by the species in $W$. Each reaction involving the species from $W^{sr}$ can not take place on $F_W$, but there exists a vector $v\in \mathscr{S}$ such that $v_{j}=0$ for every $X_j\in W$, while $v_{\tilde{j}}\neq 0$.
  \item $W^{tr}$: This subset consists of species $X_{\hat{j}}$ that are entirely restricted by the species in $W$, implying that any vector $v\in\mathscr{S}$ with $v_{j}=0$ for each $X_j\in W$ must satisfy $v_{\hat{j}}=0$. 
\end{itemize}
\begin{remark}
    Combined with the description of the Example \ref{ex:4} and Remark \ref{rek:1}, we can obtain that $W^{cv}\subset W^{tr}$ and $W^{cn}\subset W^{tf}$.
\end{remark}
The following example further illustrates three subsets described above.
\begin{example}\label{ex:4ccc}
Now let us analyze the combined system $\mathcal{M}_{\mathrm{De}}$ in Example \ref{ex:1}, where the subsystems play a vital role in a biochemical process. In this case, we assume that $A$ and $A_1$ represent the same species, and $C$ and $C_1$ represent the same species.  

The system $\mathcal{M}_{\mathrm{De}}$ is a complex delayed chemical reaction network with a 6-dimensional stoichiometric subspace, and the set of intersecting species is denoted as $\mathcal{S}_c=\{A, C, L_1\}$. By the way, all the semilocking set $\mathcal{M}_{\mathrm{De}}$ are
 \begin{equation*}
     \begin{split}
         &W_1=\{A,B,C,C_2,C_3\}; W_2=\{A,C,A_2,C_2,C_3\};\\&W_3=\{A,B,C,A_2,C_2,C_3\};W_4=\{B,C,C_2,C_3\}
         \\&W_5=\{A_2,B,C,C_2,C_3\};W_6=\{A,L,L_1,B,C,C_2,C_3\};\\&W_7=\{A,L,L_1,A_2,C,C_2,C_3\}; W_8=\mathcal{S}\\
     \end{split}
 \end{equation*}
 We can see that different with the situation in Section \ref{sec:4}, $\mathcal{S}_c\cap W\neq \emptyset$.  Considering the semilocking set $W_1$, its corresponding complement is $W_1^{c}=\{L,L_1,A_2\}$.  It is straightforward to verify that the species $L$ and $L_1$ are both in $W_1^{tf}$ as the reversible reaction  $\xymatrix{L\ar @{ -^{>}}^{~\tau_1,~k_1}@< 1pt> [r]& L_1 \ar  @{ -^{>}}^{~\tau_2,~k_2}  @< 1pt> [l]}$ can take place on boundary $F_{W_1}$. $A_2$ is in $W_1^{tr}$ as each vector $v\in \mathscr{S}$ with $v_{A_2}\neq 0$ must have $v_j\neq 0$ for some $X_j\in W$, thus the value of the concentration of $A_2$ cannot change and is a determined constant on $F_W$. Further the set $W^{sr}$ is empty. 

 But for the following simple network
\begin{equation*}
\mathcal{M}_{De:} \xymatrix{X_2+X_1\ar @{ -^{>}}@< 1pt> [r]& 2X_1 \ar  @{ -^{>}}  @< 1pt> [l]};\xymatrix{X_3+X_1\ar @{ -^{>}}@< 1pt> [r]& 2X_1\ar  @{ -^{>}}  @< 1pt> [l]}
\end{equation*}
$W=\{X_1\}$ is obviously a semilocking set of $\mathcal{M}_{\mathrm{De}}$. On the boundary $F_W$, all reactions can not take place. But there exists the vector $v=(0,1,-1)\in \mathscr{S}$. Therefore, $X_2, X_3$ are all in $W^{sr}$.
\end{example}
Based on the proof of Lemma \ref{lem:c21}, we can get the extension form of results in Case II and III. 
\begin{lemma}[Extension results of Case II and III]\label{lem:c512}
  $\mathcal{M}_{\mathrm{De}}$ is a DeCBMAS consisting of 2-dimensional sub-DeCBMASs $\mathcal{M}_{\mathrm{De}}^{(p)}, p=1,\cdots,\ell$. The boundary $F_W$ cannot have an $\omega$-limit point of any trajectory starting from a positive point if each species in $W^c$ falls in $W^{tf}$ or $W^{tr}$ and $W^{tf}$ is non-empty.
\end{lemma}

%   Thus the following special case holds.
%\begin{corollary}
%Given a DeCBMAS $\mathcal{M}_{\mathrm{De}}$ composed of several sub-DeCBMASs $\mathcal{M}_{\mathrm{De}}^{(p)}, p=1,\cdots,\ell$. The 
%The stoichiometric subspace $\pi_{W}(\mathscr{S}^{(p)})$ of the subsystem $\mathcal{M}_{\mathrm{De}}^{(p)}$ are independent of each other, where $\mathcal{S}_c$ represents the set of intersecting species.
%\end{corollary}

Thus the following result about the persistent is obviously.
\begin{theorem}
A DeCBMAS $\mathcal{M}_{\mathrm{De}}$ is a combination of $\ell$ sub-DeCBMASs $\mathcal{M}_{\mathrm{De}}^{(p)}, p=1,\cdots, \ell$. $\mathcal{M}_{\mathrm{De}}$ is persistent if for each semilocking set $W\cap \mathcal{S}_c\neq \emptyset$, $W$ satisfies the condition in Lemma \ref{lem:c531} or Lemma \ref{lem:c512}.
\end{theorem}
Also, when we shift our perspective to a conservative system and utilize its reduecd version, the following conclusions can be obtained.
\begin{theorem}\label{thm:semicon}
    $\mathcal{M}_{\mathrm{De}}$ is a DeCBMAS and $W$ is its semilocking set. Let $\mathcal{\tilde{M}}_{\mathrm{De}}$ denote the reduced delayed system on $W$. If the dimension of the stoichiometric subspace $\tilde{\mathscr{S}}$ is less than the number of species in $W$, then $F_W$ of $\mathcal{M}_{\mathrm{De}}$ cannot contain an $\omega$-limit point of some trajectory starting from a positive initial point $\psi\in\mathscr{C}_+$.
\end{theorem}

\begin{example}
Consider $\mathcal{M}_{\mathrm{De}}$ in Example \ref{ex:1} with $A=A_1$ and $C=C_1$, we analyze each semilocking set which is presented in Example \ref{ex:4ccc} individually. 

(1) $W_1$ and $W_2$: the subsets $W_1^{(1)}$ and $W_2^{(1)}$ are all $\{A\}$, representing a facet-type semilocking set of the subsystem $\mathcal{M}_{\mathrm{De}}^{(1)}$; $W_1^{(2)}=\{A,B,C\}$ and $W_2^{(2)}=\{A,C,A_2\}$ are all vertex-type semilocking sets of the subsystem $\mathcal{M}_{\mathrm{De}}^{(2)}$, while $W_1^{(3)}$ and $W_2^{(3)}$ contains $\{C,C_2,C_3\}$, denoting a vertex-type semilocking set of the subsystem $\mathcal{M}_{\mathrm{De}}^{(3)}$. Species $B\in W_1^{c}$ and $A_2\in W_2^{c}$ are all in $W^{tr}$, while species $L$ and $L_1$ are actually in the set of $W^{tf}$. Besides, the set $W^{sr}$ is empty. Lemma \ref{lem:c512} guarantees the non-existence of an $\omega$-limit point of the trajectory starting from some positive point. The same analysis can be applied into $W_6$ and $W_7$.

(2) $W_4$: $W_4^{(1)}=\emptyset$.  The subsets  $W_4^{(2)}=\{B,C\}$ and $W_4^{(3)}=\{C,C_2,C_3\}$ represent a facet-type semilocking set and a vertex-type semilocking set, respectively, in their corresponding subsystems. Species $B$ is only contained in a facet-type semilocking set $W_4^{(2)}$. Thus Theorem \ref{cor:c53} can be applied in this semilocking set.

(3) For $W_3$, $W_5$, the conditions specified in Theorem \ref{thm:semicon} are satisfied. 

In conclusion, the persistence of $\mathcal{M}_{\mathrm{De}}$ is established by considering it as a combination of three sub-DeCBMASs with coupling.
\end{example}

\section{Conclusion}\label{sec:6}
In the context of the paper, combination methods provide a way to integrate smaller sub-DeCBMASs into a larger system like the Example \ref{ex:1}, expanding the class of delayed chemical reaction systems with persistence. The combination methods have close connection with the concept of modularity in synthetic biology, which share common goals and principles. By employing modular components and combination methods, researchers can create more intricate and versatile biological systems that exhibit desired properties and behaviors. Furthermore, the intersecting species considered in this paper is fairly common occurrence as a resources for inter-module competition. It constitutes a crucial factor contributing to the failure of modularity and has garnered significant attention in synthetic biology. Therefore, the inheritable results discussed in the paper not only contribute to the analysis of complex chemical reaction systems but also provide insights into the modular design principles that underpin synthetic biology.
% In the future, we will consider continuing to analyze the possibility that the boundary equilibrium point becomes the $\omega$ limit point. Further we apply these results to the analysis of dynamical properties such as persistence. }
 
% \newtext{These results on the $\omega$-limit set of chemical reaction systems help usto better understand the dynamical properties of the chemical reactionsystem, such as stability, persistence, etc.}

 % OR

%\begin{figure}
%\begin{center}
%\epsfig{file=jcaesar,width=7cm}
%\caption{Gaius Julius Caesar, 100--44 B.C.}
\begin{ack}                               % Place acknowledgements
This work was funded by China Postdoctoral Science Foundation under Grant No. 2023M733115, the National Nature Science Foundation of China under Grant No. 12071428 and 62111530247, and the Zhejiang Provincial Natural Science Foundation of China under Grant No. LZ20A010002.
\end{ack}
\section*{Appendix: Proofs of results in paper.}
\textbf{Proof of Lemma \ref{thm:4.9}:}
Assuming that $W^{(p_1)}$ is the facet-type semilocking set of $\mathcal{M}_{\mathrm{De}}^{(p_1)}$. From the proof of the Lemma \ref{lem:c3}, we can see that $x_j^{\psi}(t)$ with $X_j\in W^{(p_1)}$ can not converges to zero independent of the value of concentration of $X_j\notin W$ as long as the initial point $\psi\in\mathscr{C}_+$. Thus the result holds.$\hfill\blacksquare$

\textbf{Proof of Lemma \ref{lem:vv}:}
If the semilocking set $W$ is the species set $\mathcal{S}$, its corresponding boundary is exactly a vertex, namely the zero vector, denoted as $\mathbbold{0}_n$.
Otherwise, suppose there exists a nontrivial semilocking set $W$ is not a vertex-type semilocking set, but each subset $W^{(p)}$ is a vertex-type semilocking set. Thus there exists a nonzero vector $v\in\mathbb{R}^{n}$ such that $v_{j}=0$ for each $S_j\in W$ and $v\in\mathscr{S}$, where $\mathscr{S}$ is the stoichiometric subspace of the system $\mathcal{M}_{\mathrm{De}}$. Consequently, there exist a set of not all zero constants $k_1,\cdots,k_{\ell}$ and a set of non-zero vectors $\{v^{(i)}, i=1,\cdots,\ell\vert v^{(i)}\in \text{span}\{{y'_{.i}}^{(p)}-y^{(p)}_{.i}\}\}$ satisfying the equation
\begin{equation}\label{eq:13}
k_1v^{(1)}+\cdots+k_\ell v^{(\ell)}+v=0.
\end{equation}
Since $v_j=0$ for all $S_j\in W$, the following equation also holds:
\begin{equation}\label{eq:14}
k_1\pi_{W}(v^{(1)})+\cdots+k_\ell \pi_W(v^{(\ell)})=0.
\end{equation}
Considering the fact that $S_c\cap W$ is empty, we conclude that either $\pi_{W}(v^{(1)}),\cdots,\pi_{W}(v^{(\ell)})$ are independent of each other, or there exists at least one zero vector among them. However, since $k_1,\cdots,k_\ell$ are not all zero in \eqref{eq:13} and \eqref{eq:14}, there must be at least one $\pi_{W}(v^{(p)})$ that is a zero vector. Moreover, the non-zero vector $v^{(p)}\in \text{span}\{{y'_{.i}}^{(p)}-y^{(p)}_{.i}\}$, which implies $\pi_{\mathcal{S}^{(p)}}(v^{(p)})\in\mathscr{S}^{(p)}$ with $v^{(p)}_j=0$ for each $S_j\in W^{(p)}$. Consequently, this is contradiction with the assumption that $F_{W^{(p)}}$ is a vertex. Thus $F_W$ is also a vertex of $\mathcal{M}^{\mathrm{De}}$ and we conclude the result.$\hfill\blacksquare$

\textbf{Proof of Lemma \ref{lem:c21}:}
Without loss of generality, let us assume that $\mathcal{M}_{\mathrm{De}}^{(p)}, p=1,\dots,\ell_1$ are subsystems whose corresponding $W^{(p)}$ are non-empty subsets. If $\ell_1$ is equal to $\ell$, the result can be directly derived from Theorem \ref{thm:c1}. However, if $\ell_1$ is less than $\ell$, the semilocking set $W$ is not a locking set because the reactions in $\mathcal{M}_{\mathrm{De}}^{(p)}$ can also occur on the boundary $F_W$. Only the boundary equilibrium $\bar{x}^*$ has the possibility of being an $\omega$-limit point of some trajectory starting from a positive point. And we denote $\cup_{p=1}^{\ell_1}\mathcal{M}_{\mathrm{De}}^{(p)}\triangleq\mathcal{M}_{\mathrm{De}}^{1}$ and $\cup_{p'=\ell_1+1}^{\ell}\mathcal{M}_{\mathrm{De}}^{(p')}\triangleq\mathcal{M}_{\mathrm{De}}^2$. 

Let us denote the semilocking set $W$ of $\mathcal{M}_{\mathrm{De}}$ as $W=\{X_1,\cdots,X_d\}$, the species $X_{j'}$ in $W^{cv}$ as $\{X_{d+1},\cdots,X_{d_1}\}$, and the species $X_{l}$ in $W^{cn}$ as $\{X_{d_1+1},\cdots,X_{n}\}$. According to the definition of $W^{cv}$, there does not exist a vector of the form $v=(\mathbbold{0}_d,v_{d+1},\cdots, v_{d_1}, v_{d_1+1},\cdots, v_n)$ with some $v_i\neq 0, i=d+1,\cdots, d_1$. This is because each species $X_{j'}$ only exists in the system $\mathcal{M}_{\mathrm{De}}^{1}$, and $W$ is a vertex-type semilocking set of $\mathcal{M}_{\mathrm{De}}^1$ from the Lemma \ref{lem:vv}. 
We assume that the boundary point $\bar{x}^*$ on the boundary $F_W$ is the $\omega$-limit point of some trajectory $x^\phi(t)$ starting from a positive initial point $\phi$. Then, by considering a small enough $\epsilon_1 > 0$, we can find that there exist two distinct time scales between the dynamics of the species $X_l \in W^{cn}$ and the species $X_{j'} \in W^{cv}$ when the trajectory $x^\phi(t)$ enters the $\epsilon_1$-neighborhood of $\bar{x}^*$. These different time scales arise due to the significant difference in the reaction rates between the reactions in $\mathcal{M}_{\mathrm{De}}^1$ and those in $\mathcal{M}_{\mathrm{De}}^2$.
  
  The different time scales can be presented in the dynamics of these species written in the following form 
  \begin{equation}\label{eq:ts}
  \begin{split}
      \epsilon_1\dot{x}_l&=\epsilon_1(f^1_l(\epsilon_1)+f^2_l)\\
      \dot{x}_{j'}&=f^1_{j'}(\epsilon_1).
      \end{split}
  \end{equation}
  As each $X_l\in W^{cn}$ can exist both in the $\mathcal{M}_{\mathrm{De}}^1$ and $\mathcal{M}_{\mathrm{De}}^2$, $f^1_l,f^2_l$ in the above equation denote the changes of the concentration of species $X_l$ causing by the reactions in $\mathcal{M}_{\mathrm{De}}^1$ and $\mathcal{M}_{\mathrm{De}}^2$ respectively. Each reaction rate in $\mathcal{M}_{\mathrm{De}}^1$ contains the concentration of $X_j\in W$ which is smaller than $\epsilon_1$. Thus we can conclude that $f^1<<f^2$ and $\dot{x}_l>>\dot{x}_{j'}$ for each $X_l$ and $X_j$. In other words, once the trajectory $x^\phi(t)$ enters into the $\epsilon_1$-neighbourhood of $\bar{x}^*$, the concentration of the species $X_l\in W^{cn}$ will reach the value of equilibrium quickly comparing to $X_{j'}$. 
  As the $\bar{x}^*$ is an $\omega$-limit point of the trajectory $x^\phi(t)$,  $\psi\in\mathscr{C}^+$ denoted as $\psi(s)=x^\phi(t+s), s\in [-\tau_{\mathrm{m}},0]$ can be found which represents the point on the trajectory that firstly enters into the $\epsilon_2$-neighbourhood of $\bar{x}^*$ and another point $\hat{\psi}(s)\in\mathscr{C}^+$ on $x^\phi(t)$ in the $\epsilon_2/2$-neighbourhood of $\phi$ where $\epsilon_2<\epsilon_1$, namely,
  \begin{equation*}
      \begin{split}
    \vert \psi_j(s)\vert<\epsilon_2,~\forall j=1,...,d,\\
     \vert \psi_{j'}(s)-\phi_{j'}\vert<\epsilon_2,~\forall j'=d+1,...,d_1.~\\
     \Vert\hat{\psi}(s)-\phi\Vert<\epsilon_2/2,\\
    \hat{\psi}_{j}(s)\leq \phi_{j}(s)/2,~\forall j=1,...,d.\\
     \end{split}
 \end{equation*}
The following equation considers the value of the difference between $V(\hat{\psi})$ and $V(\psi)$ where $V$ is the Lyapunov functional of the delayed system $\mathcal{M}_{\mathrm{De}}$. And the differential mean value theorem is also applied into $V(\hat{\psi})-V(\psi)$ 
\begin{scriptsize}
\begin{equation*}
\begin{split}
      &~~~~V(\hat{\psi})-V(\psi)\\
&=\underbrace{\sum^{d}_{j=1}\left(\ln(\tilde{\psi}_{j})-\ln(\bar{x}_{j})\right)\left(\hat{\psi}_{j}(0)-\psi_j(0)\right)}_{\alpha}\\
            & +\sum\limits^{d_1}_{j'=d+1}\left(\ln(\tilde{\psi}_{j'})-\ln(\bar{x}_{j'})\right)\left(\hat{\psi}_{j'}(0)-\psi_{j'}(0)\right)\\
    &+\underbrace{\sum\limits^{n}_{l=d_{1}+1}\left(\ln(\tilde{\psi}_{l})-\ln(\bar{x}_{l})\right)\left(\hat{\psi}_{l}(0)-\psi_{l}(0)\right)}_{\alpha_1}\\
   &+\underbrace{\sum_{i=1}^{r^1}k_i\int^{0}_{-\tau_i}\left(\ln(\tilde{\varPsi}(s))-\ln(\bar{x}^{y_{.i}})\right)\left(\hat{\psi}(s)^{y_{.i}}-\psi(s)^{y_{.i}})\right)ds}_{\beta}\\
+&\underbrace{\sum_{i'=r^1+1}^{r}k_{i'}\int^{0}_{-\tau_{i'}}\left(\ln(\bar{\varPsi}(s))-\ln(\bar{x}^{y_{.i'}})\right)\left(\hat{\psi}(s)^{y_{.i'}}-\psi(s)^{y_{.i'}})\right)}_{\beta_1}ds,
\end{split}
\end{equation*}
\end{scriptsize}
where $R_i, i=1,\cdots, r^1$ and $R_{i’},i'=r^1+1,\cdots,r$ are the reactions in $\mathcal{M}_{\mathrm{De}}^1$ and $\mathcal{M}_{\mathrm{De}}^2$ respectively.  $\tilde{\psi}_{j},~\tilde{\psi}_{j'},~\tilde{\psi}_{l},~\tilde{\varPsi}(s),\bar{\varPsi}(s)$ are constrained by $\tilde{\psi}_{j}\in (\hat{\psi}_{j}(0),\psi_{j}(0))$,$\tilde{\psi}_{j'}\in (\hat{\psi}_{j'}(0),\psi_{j'}(0))$ $\tilde{\psi}_{l}\in (\hat{\psi}_{l}(0),\psi_{l}(0))$, $\varPsi(s)\in(\hat{\psi}(s)^{y_{.i}},\psi(s)^{y_{.i}})$, and $\bar{\varPsi}(s)\in(\hat{\psi}(s)^{y_{.i}},\psi(s)^{y_{.i}})$. 
As the same proof as that in \cite{Zhang2021}, ($\alpha$), ($\beta$) are all larger than zero. And further from \eqref{eq:ts} and the quick balance of $\mathcal{W}^{cn}$, there exists $\hat{\psi}_{l}(0)-\psi_{l}(0)=o(\epsilon_1)$ is an infinitesimal of $\hat{\psi}_{j'}(0)-\psi_{j'}(0)$ and $\ln(\tilde{\psi}_{l})-\ln(\bar{x}_{l})$ is bounded. As each reaction in $\mathcal{M}_{\mathrm{De}}^2$ does not exist the species in $W$, each $\hat{\psi}(s)^{y_{.i'}}-\psi(s)^{y_{.i'}}$ is also an infinitesimal of $\hat{\psi}_{j'}(0)-\psi_{j'}(0)$ and $\ln(\bar{\Psi}(s))-\ln(\bar{x}^{y_{.i'}})$ is also bounded. Thus $\alpha_1,\beta_1$ are the infinitesimal of other terms. Therefore, the dissipation of $V$ ensures that 
\begin{equation*}
\begin{split}
    \left\vert \alpha \right\vert
    &\leq \left\vert\sum^{d_1}_{j'=d+1}\left(\ln(\tilde{\psi}_{j'}))-\ln(\bar{x}_{j'})\right)\Delta \psi^{\epsilon}_{j'}\right\vert,\\
    \left\vert \beta\right\vert &\leq \left\vert\sum^{d_1}_{j'=d+1}\left(\ln(\tilde{\psi}_{j'}))-\ln(\bar{x}_{j'})\right)\Delta \psi^{\epsilon}_{j'}\right\vert.
  \end{split}
\end{equation*}
Similar to the proof of \cite{Zhang2021}, we can define the vector sequence $\{z^1(\epsilon_n)=g(\hat{\psi}_n)-g(\psi_n)\}$ which supports 
 \begin{equation*}
     z^1_{j}(\epsilon_n)\to 0, z^1_{l}(\epsilon_n) \to 0 ~,~j\in \{1,\cdots,d\}, l\in \{d_1+1,\cdots, n\}
 \end{equation*}
when $\epsilon_{n}\rightarrow 0$. And for each $n$, the vector $z^1(\epsilon_n)$ is contained in a compact space 
\begin{equation*}
    O=\{z:D_1\leq \vert z\vert_{1}\leq D_2\}\cap \mathscr{S}.
\end{equation*}
where $D_1$ and $D_2$ are all constants. Thus there exists a subsequence $\{z^1(\epsilon_{n_p})\}$ converging to the point $z^0$ as $p\rightarrow +\infty$ where $\vert z^0\vert_1>D_1$. Thus $Z^0\in\mathscr{S}$ must be $(\mathbbold{0}_d,z^0_{d+1},\cdots,z^0_{d_1},\mathbbold{0}_{n-d_1})$ with some $z^0_{j'}>0$. This contradicts that $W$ is a vertex-type semilocking set of $\mathcal{M}_{\mathrm{De}}^1$. Thus the assumption that $\bar{x}^*$ is an $\omega$-limit point of the trajectory $x^\phi(t)$ is not true. Hence, we conclude our result.$\hfill\blacksquare$

%\bibliographystyle{plain}        % Include this if you use bibtex 
%\bibliography{autosam}           % and a bib file to produce the 
                                 % bibliography (preferred). The
                                 % correct style is generated by
                                 % Elsevier at the time of printing.
\textbf{Proof of Lemma \ref{lem:c22}:}
In order to research the ability of $F_W$ to contain an $\omega$-limit point of some trajectory with a positive initial point, only the dynamics of species $X_j$ in the semilocking set $W$ make sense. Through reducing the system on $W$, the corresponding reduced system $\mathcal{\tilde{M}}$ can be divided into several independent generalized mass action systems $\mathcal{\tilde{M}}$. Then the objective can be obtained by researching whether the origin can be one $\omega$-limit point of some trajectory with a positive initial point of the 2d reduced subsystem $\mathcal{\tilde{M}}^{(p)}$.   

The following part contributes to the case of the 2d $\mathcal{\tilde{M}}^{(p)}$ with more than $2$ species, namely, $\tilde{n}^{(p)}>2$. We claim that there must exist some non-empty and non-negative vector $a_p\in \tilde{\mathscr{S}}^{\perp}$. If not, all the vectors in the form of $(0,\cdots,0,1,0,\cdots,0)\in\mathbb{R}^{\tilde{n}^{(p)}}$ are all in $\mathscr{\tilde{S}}$. This is obviously contradiction with $\dim{\mathscr{\tilde{S}}}=2$. Further non-negative vector $a_p\in\tilde{\mathcal{S}}^{\perp}$ ensures that $a_p^{\top}g(\psi)\neq a_p^{\top}g(0)$ for each initial condition $\psi\in\mathscr{C}_+$. This means that the origin is not in the $\mathcal{D}_\psi$ from the Definition \ref{def:scc}. Thus the origin can not be the $\omega$-limit point of the trajectory $x^{\psi}(t)$ starting from $\psi\in \mathscr{C}_+$.

 So we just need to consider the case that each $\tilde{\mathcal{M}}_{\mathrm{De}}^{(p)}$ is a two species system.
 As the concentration of each species is bounded, we can derive that $\tilde{k}^{(p)}_i(t)$ in $\mathcal{M}_{\mathrm{De}}^{(p)}$ is also bounded. And for each $x$ near the boundary $F_W$, $x_j$ is also have a positive lower bounded. This means that there exist two positive constants $b^{(p)}_i$ and $B^{(p)}_i$, such that $b^{(p)}_i<\tilde{k}^{(p)}_i(t)<B^{(p)}_i$. Note that the following proof focuses on the single GDeMAS $\mathcal{\tilde{M}}^{(p)}$, thus we directly use $x$ to denote $\tilde{x}^{(p)}$ for the sake of simplifying writing.
 
 Suppose that the origin is an $\omega$-limit point of some trajectory $x^\psi(t)$ with positive initial point $\psi$. Then for each $\epsilon>0$ and $t_0>0$, there exists $t>t_0$ such that $x_j(t)<\epsilon$ for each $X_j\in \mathcal{\tilde{S}}^{(p)}$. Now we consider the dynamics of $\dot{x}$ where $x$ is in the $\epsilon$-neighbourhood of the origin. 
 \begin{equation*}
 \begin{split}
     \dot{x}=&\sum_{i=1}^r\tilde{k}_i(t-\tau_i)x^{y_{.i}}(t-\tau_i)y'_{.i}-\sum_{i=1}^r\tilde{k}_i(t)x^{y_{.i}}(t)y_{.i}\\
     =&\sum_{i=1}^r\tilde{k}_i(t)x^{y_{.i}}(t)\left(\frac{\tilde{k}_i(t-\tau_i)x^{y_{.i}}(t-\tau_i)}{\tilde{k}_i(t)x^{y_{.i}}(t)}y'_{.i}-y_{.i}\right)\\
     =&\sum_{i=1}^r\tilde{k}_i(t)x^{y_{.i}}(t)\left((1+m_i(\epsilon))y'_{.i}-y_{.i}\right)\\
       =&\sum_{i=1}^r\tilde{k}_i(t)x^{y_{.i}}(t)\left(y'_{.i}-y_{.i}+m_i(\epsilon)y'_{.i}\right)\\
     \end{split} 
    \end{equation*}
  The third equation holds as only equilibrium can be the $\omega$-limit point for each trajectory of $\mathcal{M}^{(p)}$ and each $x_j(t-\tau_i)/x_j(t)$ tends to 1 when the trajectory $x^\psi(t)$ tends to the any equilibrium. $\dot{x}$ is the linear combination of vectors $y'_{.i}-y_{.i}+m_i(\epsilon)y'_{.i}$. Thus in order to verify the sign of $\dot{x}$ when $x$ is in the neighbourhood of zero, it is necessary to know which vectors have the large enough coefficients $\tilde{k}_i(t)x^{y_{.i}}(t)$ to determine the sign of $\dot{x}$ regardless of the directions of other vectors.
  
  Now we consider the relation between $x_1(t)$ and $x_2(t)$, if $\lim\limits_{x\to 0} \frac{x_1^g(t)}{x_2(t)}=0$ for each positive constant $g$, each $x^{y_{.i}}(t)$ can be expressed by $c^{ay_{1i}+by_{2i}}$ where $c(t)=x_2(t), a=\delta<<1, b=1$; if $\lim\limits_{x\to 0} \frac{x_1(t)}{x_2^g(t)}=\infty$ for each positive constant $g$, $x^{y_{.i}}(t)$ can be expressed by $c^{ay_{1i}+by_{2i}}$ where $c(t)=x_1(t), a=1, b=\delta<<1$; Also, if there exists some $g$ such that $\lim\limits_{x\to 0} \frac{x_1^g(t)}{x_2(t)}=u$ where $u$ is a positive constant, $x^{y_{.i}}(t)$ can be expressed by $c_ic^{ay_{1i}+by_{2i}}$ where $a,b,c_i>0$ are all constants. Let $\nu$ denote the vector $(a,b)$.
   Thus we can obtain that for each $x$ in different region of the $\epsilon$-neighbourhood of the origin, its dynamics can be written as 
   \begin{equation}\label{eq:eq19}
    \dot{x}=\sum_{i=1}^r\tilde{k}^c_i(t)c(t)^{\nu y_{.i}}\left(y'_{.i}-y_{.i}+m_i(\epsilon)y'_{.i}\right)
   \end{equation}
   where $c(t)\to 0$ when $x$ tends to the origin and $\tilde{k}^c_i(t)=\tilde{k}_i(t)c_i$. We partition the complex set into several equivalent classes with each complex in $\mathcal{EC}^{f}$ sharing the same value of $\nu y_{.i}=m_f$ and $m_{f_1}<m_{f_2}$ for $f_1<f_2$. It is obviously that the smaller $f$ is, the larger coefficients of  vectors $y'_{.i}-y_{.i}+m_i(\epsilon)y'_{.i}$ with $y_{.i}\in \mathcal{EC}^f$ are. In other words, these vectors are more decisive.

   (1) If there exists a reaction $R_i:y_{.i}\to y'_{.i}$ such that $y_{.i}\in \mathcal{EC}^1$ and $y'_{.i}\notin \mathcal{EC}^1$, then $a\dot{x}_1+b\dot{x}_2$ can be written as
   \begin{equation*}
   \begin{split}
      \nu \dot{x}&=\sum_{i=1}^r\tilde{k}_i(t)x^{y_{.i}}(t)\left[\nu(y'_{.i}-y_{.i})+\nu m_i(\epsilon)y'_{.i}\right]\\
      &=\sum_{f}\sum_{y_{.i}\in \mathcal{EC}^f} \tilde{k}^c_i(t)c^{m_f}(t)\left[\nu(y'_{.i}-y_{.i})+\nu m_i(\epsilon)y'_{.i}\right]\\
\end{split}
   \end{equation*}
   As $\lim\limits_{x\to 0}\frac{c^{m_f}}{c^{m_1}}=0$ for each $m_f\neq m_1$, we can find a corresponding $\epsilon_1>0$ such that ${c(t)}^{m_1}>\delta_1 c(t)^{m_f}$ in the $\epsilon_1$-neighbourhood of the origin for each constant $\delta_1>0$. Thus $\nu \dot{x}$ shares the same sign with 
   $$\nu\dot{x}^{\mathcal{EC}^1}\triangleq\sum_{y_{.i}\in \mathcal{EC}^1} \tilde{k}^c_i(t)c^{m_1}(t)\left[\nu(y'_{.i}-y_{.i})+\nu m_i(\epsilon)y'_{.i}\right].$$ And it can be written as 
   \begin{equation*}
   \begin{split}
     & \nu\dot{x}^{\mathcal{EC}^1} =\sum_{\substack{y_{.i}\in \mathcal{EC}^1\\y'_{.i}\in \mathcal{EC}^1}} \tilde{k}^c_i(t)c^{m_1}(t)\left[\nu(y'_{.i}-y_{.i})+\nu m_i(\epsilon)y'_{.i}\right]
       \\+&\sum_{\substack{y_{.i}\in \mathcal{EC}^1\\y'_{.i}\notin \mathcal{EC}^1}} \tilde{k}^c_i(t)c^{m_1}(t)\left[\nu(y'_{.i}-y_{.i})+\nu m_i(\epsilon)y'_{.i}\right]\\
       =&c(t)^{m_1}\left[\sum_{\substack{y_{.i}\in \mathcal{EC}^1\\y'_{.i}\notin \mathcal{EC}^1}} \tilde{k}^c_i(t)\nu(y'_{.i}-y_{.i})+\sum_{y_{.i}\in \mathcal{EC}^1}\tilde{k}^c_i(t)\nu m_i(\epsilon)y'_{.i}\right]\\
       \end{split}
       \end{equation*}
       \begin{equation*}
       \begin{split}
       =&c(t)^{m_1}\left[\underbrace{\sum_{\substack{y_{.i}\in \mathcal{EC}^1\\y'_{.i}\notin \mathcal{EC}^1}} \tilde{k}^c_i(t)(m_f-m_1)}_{A}+\underbrace{\sum_{y_{.i}\in \mathcal{EC}^1}\tilde{k}_i^c(t)\nu m_i(\epsilon)y'_{.i}}_B\right]
       \end{split}
   \end{equation*}
   The first part of the equation is large than zero due to $m_f>m_1$ for each $f$. For the second part, we can find an $\epsilon_2>0$, such that $A>>\vert B\vert $. So when $\epsilon=\min\{\epsilon_1,\epsilon_2\}$, there exists $\nu \dot{x}>0$. Combining with the positiveness of the elements of $\nu$, there must exist one  species $X_j$ such that $\dot{x}_j>0$. In other words, the trajectory will go away from the origin once it enters into the $\epsilon$-neighbourhood of the origin.
   
(2) If for each reaction $R_i: y_{.i}\to y'_{.i}$ with $y_{.i}\in \mathcal{EC}^1$, $y'_{.i}\in \mathcal{EC}^1$ also holds, meaning that each reaction $y'_{.i}\to y_{.i}$ satisfies that $(y'_{.i}-y_{.i})\cdot \nu=0$. 
Thus the span of all the reaction vectors $y'_{.i}-y_{.i}$ with $y_{.i}\in \mathcal{EC}^1$ is actually the orthogonal complement space of $\text{span}{\nu}$ respect to $\mathbb{R}^2$. 
 Thus the set of all the reaction vectors of the reactions in $\mathcal{ER}^{1}$ is $1d$ and a vector $w$ is its basis such that each $y'_{.i}-y_{.i}$ with $y_{.i}\in\mathcal{EC}^{1}$ can be denoted as $p_iw$ where $p_i$ is a constant. 
Then we can further rewrite the change caused by the reactions with $y_{.i}\in\mathcal{EC}^1$ denoted as $\dot{x}^{\mathcal{EC}^1}$
\begin{equation}\label{eq:ec1}
\begin{split}
   \dot{x}^{\mathcal{EC}^1} =&\sum_{y_{.i}\in\mathcal{EC}^1}k_i(t)c_ic(t)^{\nu y_{.i}}(p_iw+m_i(\epsilon)y'_{.i})\\
    =&c(t)^{m_1}(\sum_{p_i>0}k_i(t)c_ip_i-\sum_{p_i<0}k_i(t)c_i\vert p_i\vert)w\\
    +&\sum_{y_{.i}\in\mathcal{EC}^1}^rk_i(t)c_ic(t)^{\nu y_{.i}}m_i(\epsilon)y'_{.i}
    \end{split}
    \end{equation}
As the origin system $\mathcal{M}_{\mathrm{De}}$ is a DeCBMAS, each $\omega$-limit point is an equilibrium. Thus from the definition of the $k_i(t)$, we can obtain that $k_i(t)$ equal to some positive constant $\kappa_i$ when $x(t)$ is the origin of the reduced system.

If $\sum_{p_i>0}\kappa_ic_ip_i w\neq \sum_{p_i<0}\kappa_ic_i\vert p_i\vert w$ and $\sum_{p_i>0}\kappa_ic_ip_i -\sum_{p_i<0}\kappa_ic_ip_i=M$, the following equation holds
\begin{equation}
    \sum_{p_i>0}k_i(t)c_ip_i -\sum_{p_i<0}k_i(t)c_ip_i=M'\neq 0
\end{equation}
when $x(t)$ is in the $\epsilon_3$-neighbourhood of the origin for a small enough constant $\epsilon_3>0$. Then the first part of the $\dot{x}^{\mathcal{EC}^1}$ is not equal to zero in this case. %Further noting that $w$ cannot in a form of $w=(w_1,0)$ or $w=(0,w_2)$ or $w=(0,0)$, otherwise $y_i$ and $y'_{.i}$ can not in the $\mathcal{EC}^1$.
From the fact that $w\in \text{span}{\nu}^{\perp}$, we can directly conclude that $w$ is in the form of $w=(w_1,w_2)$ where $w_1w_2<0$. 
So an $\epsilon_4$ can be found such that the absolute value of the second term in Equation \ref{eq:ec1} is much smaller than that of the first term when $x(t)$ is in the $\epsilon_4$-neighbourhood of the origin. Combining with the form of the $w$, there must exist some species $X_j$ such that $\dot{x}^{\mathcal{EC}^1}_j(t)>0$ when $x(t)$ is in the $\epsilon_4$-neighbourhood of the origin. Further choosing $\epsilon=\min\{\epsilon_1, \epsilon_3,\epsilon_4\}$, we can conclude that $\dot{x}_j>0$ when $x(t)$ enters into the $\epsilon$-neighbourhood of the origin. This contradicts to the assumption that the origin is the $\omega$-limit point of some trajectory starting from a positive point.

The following part contributes to the proof of the case that  $\sum_{p_i>0}\kappa_ic_ip_i w = \sum_{p_i<0}\kappa_ic_ip_i w$, namely, the first part of $\dot{x}^{\mathcal{EC}^1}$ equal to zero. 
In this case, the sign of $\dot{x}^{\mathcal{EC}^1}$ can not reflect the sign of $\dot{x}$ of the reduced system $\mathcal{\tilde{M}}^{(p)}$. For each $\nu$, as the $\text{span}\{\nu\}^{\perp}$ is 1d, there must there exist reactions such that $y_{.i}$ and $y'_{.i}$ are not in the same $\mathcal{EC}^f$. And we denote $f_1$ as the minimal $f$ such that a reaction $y_{.i}\to y'_{.i}$ with $y_{.i}\in \mathcal{EC}^{f_1}, y'_{.i}\in\mathcal{EC}^{f_2}, f_1<f_2$ can be found. Then we can find a corresponding $\epsilon_5>0$ such that the direction of the trajectory of $\dot{x}$ can be determined by the direction of $\sum_{f=1}^{f_1}\dot{x}^{\mathcal{EC}^f}$ when $x(t)$ is in the $\epsilon_5$-neighbourhood of the origin. And the sign of $\dot{x}^{\mathcal{EC}^{f_1}}$ can be determined by the part (1) of this proof. Also, each  $\dot{x}^{\mathcal{EC}^{f}}$ can be denoted into two parts as \eqref{eq:ec1}. The first part is zero or at least exists some species $X_j$ with $\dot{x}_j(t)>0$. And each term of the second part $\sum_{y_{.i}\in\mathcal{EC}^f}k_i(t)c_ic(t)^{\nu y_{.i}}m_i(\epsilon)y'_{.i}$ is actually in the form of 
\begin{equation}
\sum_{y_{.i}\in \mathcal{EC}^f}k_i(t-\tau_i)c_ic(t-\tau_i)^{m_f}y'_{.i}-\sum_{y_{.i}\in \mathcal{EC}^f}k_i(t)c_ic(t)^{m_f}y'_{.i}.
\end{equation}
$x(t-\tau_i)^{y_{.i}}$ can be written as $c_ic(t-\tau_i)^{m_f}$ in the above equation due to the fact that $\lim_{x\to 0}\frac{x_j(t-\tau_i)}{x_j(t)}=1$.  
 Also, as each species $X_j\notin W$ must exist in some reactions which do not contain the species in $W$ due to $W^{cv}=\emptyset$. Thus the speed of the convergence of the species $X_j\notin W$ is much larger than that of the species in $W$. Thus sign of the above equation is equivalent to 
 \begin{equation}
\sum_{y_{.i}\in \mathcal{EC}^f}\bar{k}_ic_ic(t-\tau_i)^{m_f}y'_{.i}-\sum_{y_{.i}\in \mathcal{EC}^f}\bar{k}_ic_ic(t)^{m_f}y'_{.i}.
\end{equation}
where $\bar{k}_i=k_i\prod_{j\notin W} (\bar{x}^*_j)^{y_{ji}}$ when $x(t)$ comes into $\epsilon$-neighbourhood of $\bar{x}^*$ for some small enough $\epsilon_6$ where $\bar{x}^*$ is the unique boundary equilibrium of the boundary $F_W$ of the system $\mathcal{M}_{De}$. Then let $t_0$ denote the first time $c(t_0)<\epsilon_7$ for some $\epsilon_7<\epsilon_6$, there exists $c(t_0-\tau_i)>c(t_0)$ for each $\tau_i$. At this time, the second term of $\dot{x}^{\mathcal{EC}^1}$ is positive for each $X_j\in W$. Thus $\dot{c}(t)>0$, the trajectory cannot enter into the smaller neighbourhood of the boundary point $\bar{x}^*$. Further combining the first term in $\dot{x}^{\mathcal{EC}^{1}}$ and $\dot{x}^{\mathcal{EC}^{f_1}}$, we conclude that the boundary $F_W$ of the $\mathcal{M}_{\mathrm{De}}$ can not exists any $\omega$-limit point of the trajectory starting from some positive initial point.  $\hfill\blacksquare$

\textbf{Proof of Theorem \ref{cor:c53}:} We also consider the reduced system of $\tilde{\mathcal{M}}_{\mathrm{De}}$ respective to $W$. The existence of some vertex-type semilocking sets does not affect the dynamics of $X_j$ as $\dot{x}_j$ is also the sum of $\dot{x}_j^{(p)}$ whose corresponding $W^{(p)}$ is a facet type semilocking set of $\mathcal{M}^{(p)}$. Thus the concentration of $X_j\in W$ can not tend to zero in this case from the proof of Lemma \ref{lem:c531}. Thus the result is obtained. $\hfill\blacksquare$

\textbf{Proof of Lemma \ref{lem:c512}:}
All the conditions ensures that the proof of the Lemma \ref{lem:c21} also holds. Thus this result can be concluded. $\hfill\blacksquare$ 

\textbf{Proof of Theorem \ref{thm:semicon}:} The reduced system of $\mathcal{M}_{\mathrm{De}}$ denoted as $\tilde{\mathcal{M}}_{\mathrm{De}}$ is conservative from the condition. Similar to the proof of $2d$ case at the beginning of Lemma \ref{lem:c22},  there exists a non-negative vector $v\in \mathbb{R}^n_{\geq 0}$ with $v_j=0$ for each $X_j\in W$. And the conservation between the species in semilocking set $W$ ensures that $F_W$ can not possess some $\omega$-limit point of trajectory with positive initial point $\psi\in \mathcal{C}^+$.

\begin{wrapfigure}{l}{0mm}
\includegraphics[width=0.95in,height=1.25in,clip,keepaspectratio]{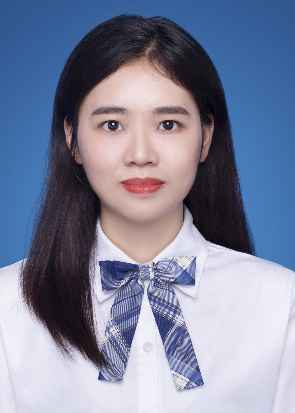}
\end{wrapfigure}
Xiaoyu Zhang received the B.Sc. degree in mathematics from Anhui University, China, in 2017, and the Ph.D. degree from Zhejiang University, China, in 2022. She is currently a postdoctor in the Department of Control Science and Engineering at Zhejiang University.
Her research interests are in the areas of control theory and chemical reaction networks theory.

\begin{wrapfigure}{l}{0mm}
     \includegraphics[width=1in,height=1.25in,clip,keepaspectratio]{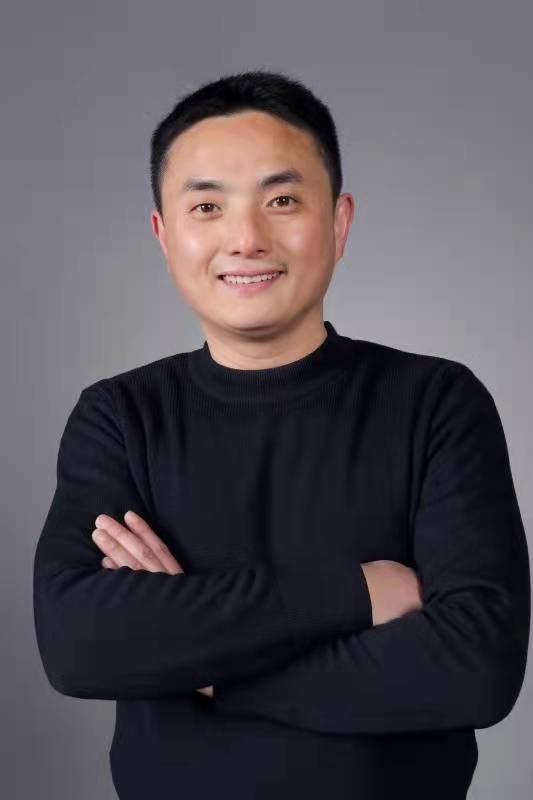}
\end{wrapfigure}
 Chuanhou Gao received the B.Sc. degrees in Chemical Engineering from Zhejiang University of Technology, China, in 1998, and the Ph.D. degrees in Operational Research and Cybernetics from Zhejiang University, China, in 2004. From June 2004 until May 2006, he was a Postdoctor in the Department of Control Science and Engineering at Zhejiang University. Since June 2006, he has joined the Department of Mathematics at Zhejiang University, where he is currently a Professor. He was a visiting scholar at Carnegie Mellon University from Oct. 2011 to Oct. 2012. 
 His research interests are in the areas of data-driven modeling, control and optimization, chemical reaction network theory and thermodynamic process control. He is an associate editor of IEEE Transactions on Automatic Control and of International Journal of Adaptive Control and Signal Processing.

\begin{wrapfigure}{l}{0mm}
	\includegraphics[width=1in,height=1.25in,clip,keepaspectratio]{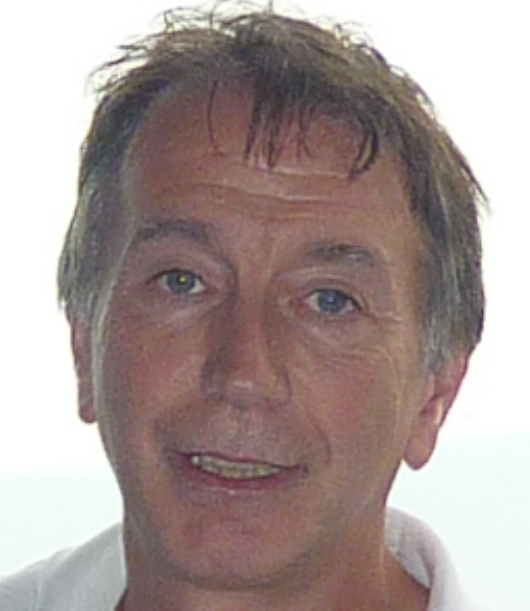}
\end{wrapfigure}
Denis Dochain received his degree in Electrical engineering in 1982 from the Université Catholique de Louvain, Belgium. He completed his Ph.D. thesis and a ``th\`{e}se d'agr\'{e}gation de l'enseignement sup\'{e}rieur" in 1986 and 1994, respectively, also at the Université Catholique de Louvain, Belgium. He has been CNRS associate researcher at the LAAS (Toulouse, France) in 1989, and Professor at the Ecole Polytechnique de Montr\'{e}al, Canada in 1987-88 and 1990-92. He has been with the FNRS (Fonds National de la Recherche Scientifique, National Fund for Scientific Research), Belgium since 1990. Since September 1999, he is Professor at the ICTEAM (Institute), Universit\'{e} Catholique de Louvain, Belgium, and Honorary Research Director of the FNRS. He has been invited professor at Queen's University, Kingston, Canada between 2002 and 2004. He is full professor at the UCL since 2005. He is the Editor-in-Chief of the Journal of Process Control, senior editor of the IEEE Transactions of Automatic Control and associate editor of Automatica. He is active in IFAC since 1999 (Council member, Technical Board member, Publication Committee member and chair, TC and CC chair). He received the IFAC outstanding service award in 2008 and is an IFAC fellow since 2010. He received the title of Doctor Honoris Causa from the INP Grenoble on December 13, 2020.
His main research interests are in the field of nonlinear systems, thermodynamics based control, parameter and state estimation, adaptive extremum seeking control and distributed parameter systems, with application to microbial ecology, environmental, biological and chemical systems, and electrical and mechanical systems. He is the (co-)author of 5 books, more than 160 papers in refereed journals and 260 international conference papers.
\end{document}